\documentclass{article}
\usepackage{amsfonts}
\usepackage{amsmath}
\usepackage{amssymb}
\usepackage{amsthm}
\usepackage[margin=1.0in]{geometry}
\usepackage[utf8]{inputenc}
\usepackage[dvipsnames]{xcolor}
\usepackage{marginnote}
\usepackage{mathrsfs}

\newtheoremstyle{custom}
{6pt}
{3pt}
{\sl}
{}
{\bf}
{.}
{.5em}
{}
\theoremstyle{custom}
\newtheorem{theorem}{Theorem}[section]
\newtheorem{conjecture}[theorem]{Conjecture}

\title{Counterexamples to Siegel's Conjecture}
\author{Kean P.~Fallon, Madisyn Janusiak, Edward D.~Kim, Avery McLain\\\sl{University of Wisconsin-La Crosse}}
\date{}

\begin{document}
\maketitle

\begin{abstract}
We prove that the intersection of a Hirsch polytope and a cube may be a non-Hirsch polytope.
\end{abstract}

\section{Introduction}

The Hirsch Conjecture states that for any $d$-dimensional polytope $P$ with $n$ facets the diameter of the graph of $P$ is less than or equal to $n-d$. This conjecture was posed in a letter by Warren M. Hirsch to George B. Dantzig in 1957. It was motivated by linear programming since the diameter of a polytope provides a lower bound on the number of steps needed by the simplex method. The Hirsch Conjecture is known to hold in a few cases: the conjecture holds when $d \leq 3$ due to a theorem of Klee (see~\cite{Klee:PathsII}) in 1970, dual transportation polytopes (see~\cite{Balinski:DualTransportation}), $0$-$1$ polytopes (see~\cite{Naddef:Hirsch01}), classical transportation polytopes (see~\cite{Borgwardt:network-flow-polytopes}), and when $n-d \leq 6$ (see~\cite{BremnerSchewe:nd6}), among others.

In 2010, Santos (see~\cite{Santos:CounterexampleHirsch}) found a counterexample to the conjecture with a $43$-polytope containing $86$ facets whose diameter is at least $44$. In 2011, Matschke, Santos, and Weibel found a $20$-dimensional counterexample with $40$ facets whose diameter is exactly $21$. Their $20$-dimensional counterexample to the Hirsch Conjecture is a so-called spindle obtained by taking the polar of a certain $20$-dimensional prismatoid. 

In light of the counterexamples to the Hirsch Conjecture, it would be interesting to revisit other conjectures referencing the Hirsch Conjecture. Question 17 on page 398 of~\cite{Yemelichev:Polytopes} states Siegel's Conjecture to be the following: ``If the maximum diameter conjecture is true for a given polytope then it is true for its intersection with a cube.'' Restated, the conjecture is:
\begin{conjecture}[Siegel's Conjecture]\label{SConj}
If $H$ satisfies the Hirsch Conjecture and $C$ is a cube, then $H \cap C$ satisfies the Hirsch Conjecture.
\end{conjecture}
We could not find any other reference to Siegel's Conjecture, but believe that it is motivated by restricting a linear program to have the value of all decision variables between $0$ and $1$. This article constructs counterexamples to Siegel's Conjecture.

\section{Background and notation}

A polyhedron $P$ is the set of points satisfying a finite collection of non-strict linear inequalities, and its dimension $d$ is the dimension of its affine hull. If $P$ is a subset of $\{x : \langle a,x \rangle \leq b\}$, then $P \cap \{x : \langle a,x \rangle = b\}$ is a face of $P$. The faces of $P$ of dimensions $0$, $1$, and $d-1$ (respectively) are called the vertices, edges, and facets (respectively) of $P$. Given a polyhedron $P$, there are two combinatorial objects arising from $P$ that are of interest in this article.

First, the vertices and edges of a polyhedron $P$ form a finite graph $G(P)$. The maximal distance between any pair of vertices in $G(P)$ is the diameter of the graph. The Hirsch Conjecture asserted that this diameter is at most the difference between the number of facets of $P$ and its dimension. We will refer to a polyhedron as Hirsch if it satisfies the Hirsch Conjecture and non-Hirsch otherwise.

Second, a much finer object is the face lattice of a polyhedron. Given a polyhedron $P$, the set of all faces of $P$ forms a poset $L(P)$ called the face lattice of $P$, with elements of $L(P)$ ordered by inclusion. We say two polyhedra $C$ and $D$ are combinatorially equivalent if their face lattices are isomorphic.

A bounded polyhedron is called a polytope. While the Hirsch Conjecture originally concerned polyhedra in general, in the last several decades, the conjecture has only been considered for bounded polytopes. For more on polyhedra, please see~\cite{Ziegler:Lectures}.

Let $N$ be a polar of the wide prismatoid with $20$ vertices defined in \cite{MatschkeSantosWeibel:Width5Prismatoids}. That is, $N$ can be described as the non-Hirsch polytope in $\mathbb{R}^{20}$ defined by the inequalities $\mathcal{I} = \{I_j : j \in \{0,\dots,39\}\}$, where $I_j$ is the inequality $1+ v_j{}^T x \geq 0$ and
{\footnotesize
\begin{align*}
v_0 &= (100, 0, 0, 21, -7, 0, 0, 0, 0, 0, 0, 0, 0, 0, 0, 1, 1, 0, 0, 0), \\
v_1 &= (100,0,0,16,-15,0,0,0,0,0,0,0,0,0,0,0,1,1,0,0), \\
v_2 &= (100,0,0,0,-32,0,0,0,0,0,0,0,0,0,0,0,0,1,1,0), \\
v_3 &= (100,0,0,-16,-15,0,0,0,0,0,0,0,0,0,0,0,0,0,1,1), \\
v_4 &= (100,0,0,-21,-7,0,0,0,0,0,0,0,0,0,0,0,0,0,0,0), \\
v_5 &= (100,0,0,-20,-4,0,0,0,0,0,0,0,0,0,0,0,0,0,0,1), \\
v_6 &= (100,0,0,0,32,0,0,0,0,0,0,0,0,0,0,0,0,0,0,0), \\
v_7 &= (100,0,0,20,-4,0,0,0,0,0,0,0,0,0,0,0,0,0,0,0), \\
v_8 &= (100,3/100,-1/50,0,-30,0,0,0,0,0,0,0,0,0,0,0,0,0,0,0), \\
v_9 &= (100,-3/100,-1/50,0,30,0,0,0,0,0,0,0,0,0,0,0,0,0,0,0), \\
v_{10} &= (100,-3/2000,7/2000,0,318/10,0,0,0,0,0,0,0,0,0,1,0,0,0,0,0), \\
v_{11} &= (100,3/2000,7/2000,0,-318/10,10^7,10^7,10^7,10^{10},10^{11},10^{11},10^{11},10^{11},1,0,0,0,0,0,0), \\
v_{12} &= (100,3/2000,7/2000,0,-318/10,-10^7,0,0,0,0,0,0,0,1,0,0,0,0,0,0), \\
v_{13} &= (100,3/2000,7/2000,0,-318/10,10^7,-10^7,0,0,0,0,0,0,1,0,0,0,0,0,0), \\
v_{14} &= (100,3/2000,7/2000,0,-318/10,10^7,10^7,-10^7,0,0,0,0,0,1,0,0,0,0,0,0), \\
v_{15} &= (100,3/2000,7/2000,0,-318/10,10^7,10^7,10^7,-10^{10},0,0,0,0,1,0,0,0,0,0,0), \\
v_{16} &= (100,3/2000,7/2000,0,-318/10,10^7,10^7,10^7,10^{10},-10^{11},0,0,0,1,0,0,0,0,0,0), \\
v_{17} &= (100,3/2000,7/2000,0,-318/10,10^7,10^7,10^7,10^{10},10^{11},-10^{11},0,0,1,0,0,0,0,0,0), \\
v_{18} &= (100,3/2000,7/2000,0,-318/10,10^7,10^7,10^7,10^{10},10^{11},10^{11},-10^{11},0,1,0,0,0,0,0,0), \\
v_{19} &= (100,3/2000,7/2000,0,-318/10,10^{7},10^{7},10^{7},10^{10},10^{11},10^{11},10^{11},-10^{11},1,0,0,0,0,0,0), \\
v_{20} &= (-100,30,0,0,0,0,0,1,1,0,0,0,0,0,0,0,0,0,0,0), \\
v_{21} &= (-100,4,-15,0,0,0,0,0,1,1,0,0,0,0,0,0,0,0,0,0), \\
v_{22} &= (-100,0,-33/2,0,0,0,0,0,0,1,1,0,0,0,0,0,0,0,0,0), \\
v_{23} &= (-100,-1,-16,0,0,0,0,0,0,0,1,1,0,0,0,0,0,0,0,0), \\
v_{24} &= (-100,-55/2,0,0,0,0,0,0,0,0,0,1,1,0,0,0,0,0,0,0), \\
v_{25} &= (-100,-17,18,0,0,0,0,0,0,0,0,0,1,0,0,0,0,0,0,0), \\
v_{26} &= (-100,0,38,0,0,0,0,0,0,0,0,0,0,0,0,0,0,0,0,0), \\
v_{27} &= (-100,22,17,0,0,0,0,0,0,0,0,0,0,0,0,0,0,0,0,0),
\\
v_{28} &= (-100,-10,0,1/5,-1/5,0,0,0,0,0,0,0,0,0,0,0,0,0,0,0), \\
v_{29} &= (-100,2999/100,0,-3/25,-1/5,0,0,1,0,0,0,0,0,0,0,0,0,0,0,0), \\
v_{30} &= (-100,299999/10000,0,0,1/100,0,1,0,0,0,0,0,0,0,0,0,0,0,0,0), \\
v_{31} &= (-100,-2745/100,0,1/5000,1/800,1,0,0,0,0,0,0,0,0,0,0,0,0,0,0), \\
v_{32} &= (-100,-27,0,1/500,-1/80,0,0,0,0,0,0,0,0,10^5,10^7,10^7,10^7,10^8,10^8,10^9), \\
v_{33} &= (-100,-27,0,1/500,-1/80,0,0,0,0,0,0,0,0,-10^5,0,0,0,0,0,0), \\
v_{34} &= (-100,-27,0,1/500,-1/80,0,0,0,0,0,0,0,0,10^5,-10^7,0,0,0,0,0), \\
v_{35} &= (-100,-27,0,1/500,-1/80,0,0,0,0,0,0,0,0,10^5,10^7,-10^7,0,0,0,0), \\
v_{36} &= (-100,-27,0,1/500,-1/80,0,0,0,0,0,0,0,0,10^5,10^7,10^7,-10^7,0,0,0), \\
v_{37} &= (-100,-27,0,1/500,-1/80,0,0,0,0,0,0,0,0,10^5,10^7,10^7,10^7,-10^8,0,0), \\
v_{38} &= (-100,-27,0,1/500,-1/80,0,0,0,0,0,0,0,0,10^5,10^7,10^7,10^7,10^8,-10^8,0), \\
v_{39} &= (-100,-27,0,1/500,-1/80,0,0,0,0,0,0,0,0,10^5,10^7,10^7,10^7,10^8,10^8,-10^9). \\
\end{align*}}
Starting with $N$, the polyhedron resulting from removing exactly one of the inequalities $I_j$ in $\mathcal{I}$ is a bounded polytope if and only if \[j \in \{0,1,2,3,4,5,6,7,8,9,10,11,20,21,22,23,24,25,26,27,28,29,30,31,32\}.\]
A combinatorial $d$-cube is a polytope that is combinatorially equivalent to the axis-parallel $d$-cube with vertex set $\{0,1\}^d$. To construct each counterexample to Conjecture~\ref{SConj}, we will construct a Hirsch polytope $H$ and a combinatorial cube $C$ such that $H \cap C$ is the polytope $N$ described above. By applying a projective transformation, we may assume that $C$ is a combinatorial cube in constructions for Conjecture~\ref{SConj}.

\section{Our counterexamples}
Our counterexamples to Conjecture~\ref{SConj} are obtained by constructing a cube $C$ and selecting a subset of $39$ of the inequalities in $\mathcal{I}$ to define a Hirsch polytope $H$ so that $N = H\cap C$.

\subsection{Counterexample based on inequality $0$}
If $H$ is the polyhedron obtained from $N$ by removing inequality $I_{0}$, then $H$ is a bounded Hirsch polytope and $N = H \cap C$, where $C$ is the combinatorial $20$-cube defined by
{\footnotesize\begin{align*}
-\frac{2361185124158778945083}{2361183241434822606848} &\leq \langle \left(0,\,0,\,0,\,0,\,0,\,0,\,1,\,0,\,0,\,0,\,0,\,0,\,0,\,0,\,0,\,0,\,0,\,0,\,0,\,0\right) ,x \rangle \leq \frac{2361185124159496218833}{2361183241434822606848} \\
-\frac{40066588727974991}{36028797018963968} &\leq \langle \left(0,\,0,\,0,\,1,\,0,\,0,\,0,\,0,\,0,\,0,\,0,\,0,\,0,\,0,\,0,\,0,\,-21,\,0,\,0,\,0\right) ,x \rangle \leq \frac{40081012197285329}{36028797018963968} \\
-\frac{10805905420429795}{2251799813685248} &\leq \langle \left(-100,\,0,\,0,\,-21,\,7,\,0,\,0,\,0,\,0,\,0,\,0,\,0,\,0,\,0,\,0,\,-1,\,-1,\,0,\,0,\,0\right) ,x \rangle \leq 1 \\
-\frac{4722372222305608191517}{4722366482869645213696} &\leq \langle \left(0,\,0,\,0,\,0,\,0,\,0,\,0,\,0,\,0,\,0,\,0,\,0,\,0,\,0,\,0,\,0,\,0,\,0,\,1,\,0\right) ,x \rangle \leq \frac{4722372222305608191517}{4722366482869645213696} \\
-\frac{153120913883298105}{144115188075855872} &\leq \langle \left(0,\,0,\,0,\,0,\,1,\,0,\,0,\,0,\,0,\,0,\,0,\,0,\,0,\,0,\,0,\,0,\,7,\,0,\,0,\,0\right) ,x \rangle \leq \frac{76542601517243887}{72057594037927936} \\
-\frac{18889470523027351237041}{18889465931478580854784} &\leq \langle \left(0,\,0,\,0,\,0,\,0,\,0,\,0,\,0,\,0,\,0,\,0,\,0,\,0,\,0,\,0,\,0,\,0,\,0,\,0,\,1\right) ,x \rangle \leq \frac{18889470523027351237041}{18889465931478580854784} \\
-\frac{9444738705175253405213}{9444732965739290427392} &\leq \langle \left(0,\,0,\,0,\,0,\,0,\,0,\,0,\,0,\,0,\,0,\,0,\,0,\,0,\,0,\,0,\,0,\,0,\,1,\,0,\,0\right) ,x \rangle \leq \frac{9444738705175253405213}{9444732965739290427392} \\
-\frac{39614081257140112583132697033}{39614081257132168796771975168} &\leq \langle \left(0,\,0,\,0,\,0,\,0,\,0,\,0,\,0,\,0,\,0,\,1,\,0,\,0,\,0,\,0,\,0,\,0,\,0,\,0,\,0\right) ,x \rangle \leq \frac{9903520314285028145783175131}{9903520314283042199192993792} \\
-\frac{2417851645399615732725419}{2417851639229258349412352} &\leq \langle \left(0,\,0,\,0,\,0,\,0,\,0,\,0,\,0,\,0,\,0,\,0,\,0,\,1,\,0,\,0,\,0,\,0,\,0,\,0,\,0\right) ,x \rangle \leq \frac{2417851645399615732725419}{2417851639229258349412352} \\
-\frac{147579558414782802697}{147573952589676412928} &\leq \langle \left(0,\,0,\,0,\,0,\,0,\,0,\,0,\,0,\,0,\,0,\,0,\,0,\,0,\,1,\,0,\,0,\,0,\,0,\,0,\,0\right) ,x \rangle \leq \frac{147579558414782802697}{147573952589676412928} \\
-\frac{9671406563106733594076309}{9671406556917033397649408} &\leq \langle \left(0,\,0,\,0,\,0,\,0,\,0,\,0,\,0,\,0,\,0,\,0,\,1,\,0,\,0,\,0,\,0,\,0,\,0,\,0,\,0\right) ,x \rangle \leq \frac{9671406563106733594076309}{9671406556917033397649408} \\
-\frac{309485009827534768921207957}{309485009821345068724781056} &\leq \langle \left(0,\,0,\,0,\,0,\,0,\,0,\,0,\,0,\,0,\,1,\,0,\,0,\,0,\,0,\,0,\,0,\,0,\,0,\,0,\,0\right) ,x \rangle \leq \frac{309485009827534768921207957}{309485009821345068724781056} \\
-\frac{1180593414291149733993}{1180591620717411303424} &\leq \langle \left(0,\,0,\,0,\,0,\,0,\,0,\,0,\,0,\,0,\,0,\,0,\,0,\,0,\,0,\,0,\,1,\,-1,\,0,\,0,\,0\right) ,x \rangle \leq \frac{1180597001438626595131}{1180591620717411303424} \\
-\frac{76859732694334399}{72057594037927936} &\leq \langle \left(0,\,1,\,0,\,0,\,0,\,0,\,0,\,0,\,0,\,0,\,0,\,0,\,0,\,0,\,0,\,0,\,0,\,0,\,0,\,0\right) ,x \rangle \leq \frac{77296343773163779}{72057594037927936} \\
-\frac{1180593503442084915409}{1180591620717411303424} &\leq \langle \left(0,\,0,\,0,\,0,\,0,\,0,\,0,\,1,\,0,\,0,\,0,\,0,\,0,\,0,\,0,\,0,\,0,\,0,\,0,\,0\right) ,x \rangle \leq \frac{1180593503442084915409}{1180591620717411303424} \\
-\frac{151699951751958945}{144115188075855872} &\leq \langle \left(0,\,0,\,1,\,0,\,0,\,0,\,0,\,0,\,0,\,0,\,0,\,0,\,0,\,0,\,0,\,0,\,0,\,0,\,0,\,0\right) ,x \rangle \leq \frac{80781273703537389}{72057594037927936} \\
-\frac{2417851643085731713983019}{2417851639229258349412352} &\leq \langle \left(0,\,0,\,0,\,0,\,0,\,0,\,0,\,0,\,1,\,0,\,0,\,0,\,0,\,0,\,0,\,0,\,0,\,0,\,0,\,0\right) ,x \rangle \leq \frac{2417851643085731713983019}{2417851639229258349412352} \\
-\frac{2361185035008561037417}{2361183241434822606848} &\leq \langle \left(0,\,0,\,0,\,0,\,0,\,0,\,0,\,0,\,0,\,0,\,0,\,0,\,0,\,0,\,1,\,0,\,0,\,0,\,0,\,0\right) ,x \rangle \leq \frac{2361185035008561037417}{2361183241434822606848} \\
-\frac{291198229793782281}{288230376151711744} &\leq \langle \left(1,\,0,\,0,\,0,\,0,\,0,\,0,\,0,\,0,\,0,\,0,\,0,\,0,\,0,\,0,\,0,\,-100,\,0,\,0,\,0\right) ,x \rangle \leq \frac{582225359826457723}{576460752303423488} \\
-\frac{18889473462374272866723}{18889465931478580854784} &\leq \langle \left(0,\,0,\,0,\,0,\,0,\,1,\,0,\,0,\,0,\,0,\,0,\,0,\,0,\,0,\,0,\,0,\,0,\,0,\,0,\,0\right) ,x \rangle \leq \frac{18889473462375840754951}{18889465931478580854784} \\
\end{align*}}
\subsection{Counterexample based on inequality $1$}
If $H$ is the polyhedron obtained from $N$ by removing inequality $I_{1}$, then $H$ is a bounded Hirsch polytope and $N = H \cap C$, where $C$ is the combinatorial $20$-cube defined by
{\footnotesize\begin{align*}
-\frac{590297603932444082281}{590295810358705651712} &\leq \langle \left(0,\,0,\,0,\,0,\,0,\,0,\,0,\,0,\,0,\,0,\,0,\,0,\,0,\,0,\,0,\,0,\,1,\,-1,\,0,\,0\right) ,x \rangle \leq \frac{2361191850588767073579}{2361183241434822606848} \\
-\frac{1180593414291149733993}{1180591620717411303424} &\leq \langle \left(0,\,0,\,0,\,0,\,0,\,0,\,0,\,0,\,0,\,0,\,0,\,0,\,0,\,0,\,0,\,1,\,0,\,0,\,0,\,0\right) ,x \rangle \leq \frac{1180593414291149733993}{1180591620717411303424} \\
-\frac{9671406563106733594076309}{9671406556917033397649408} &\leq \langle \left(0,\,0,\,0,\,0,\,0,\,0,\,0,\,0,\,0,\,0,\,0,\,1,\,0,\,0,\,0,\,0,\,0,\,0,\,0,\,0\right) ,x \rangle \leq \frac{9671406563106733594076309}{9671406556917033397649408} \\
-\frac{2361185124158778945083}{2361183241434822606848} &\leq \langle \left(0,\,0,\,0,\,0,\,0,\,0,\,1,\,0,\,0,\,0,\,0,\,0,\,0,\,0,\,0,\,0,\,0,\,0,\,0,\,0\right) ,x \rangle \leq \frac{2361185124159496218833}{2361183241434822606848} \\
-\frac{4722372222305608191517}{4722366482869645213696} &\leq \langle \left(0,\,0,\,0,\,0,\,0,\,0,\,0,\,0,\,0,\,0,\,0,\,0,\,0,\,0,\,0,\,0,\,0,\,0,\,1,\,0\right) ,x \rangle \leq \frac{4722372222305608191517}{4722366482869645213696} \\
-\frac{76859732694334399}{72057594037927936} &\leq \langle \left(0,\,1,\,0,\,0,\,0,\,0,\,0,\,0,\,0,\,0,\,0,\,0,\,0,\,0,\,0,\,0,\,0,\,0,\,0,\,0\right) ,x \rangle \leq \frac{77296343773163779}{72057594037927936} \\
-\frac{18889470523027351237041}{18889465931478580854784} &\leq \langle \left(0,\,0,\,0,\,0,\,0,\,0,\,0,\,0,\,0,\,0,\,0,\,0,\,0,\,0,\,0,\,0,\,0,\,0,\,0,\,1\right) ,x \rangle \leq \frac{18889470523027351237041}{18889465931478580854784} \\
-\frac{38280228348601461}{36028797018963968} &\leq \langle \left(0,\,0,\,0,\,0,\,1,\,0,\,0,\,0,\,0,\,0,\,0,\,0,\,0,\,0,\,0,\,0,\,0,\,15,\,0,\,0\right) ,x \rangle \leq \frac{76542601517243887}{72057594037927936} \\
-\frac{39614081257140112583132697033}{39614081257132168796771975168} &\leq \langle \left(0,\,0,\,0,\,0,\,0,\,0,\,0,\,0,\,0,\,0,\,1,\,0,\,0,\,0,\,0,\,0,\,0,\,0,\,0,\,0\right) ,x \rangle \leq \frac{9903520314285028145783175131}{9903520314283042199192993792} \\
-\frac{2417851645399615732725419}{2417851639229258349412352} &\leq \langle \left(0,\,0,\,0,\,0,\,0,\,0,\,0,\,0,\,0,\,0,\,0,\,0,\,1,\,0,\,0,\,0,\,0,\,0,\,0,\,0\right) ,x \rangle \leq \frac{2417851645399615732725419}{2417851639229258349412352} \\
-\frac{147579558414782802697}{147573952589676412928} &\leq \langle \left(0,\,0,\,0,\,0,\,0,\,0,\,0,\,0,\,0,\,0,\,0,\,0,\,0,\,1,\,0,\,0,\,0,\,0,\,0,\,0\right) ,x \rangle \leq \frac{147579558414782802697}{147573952589676412928} \\
-\frac{291128156979952409}{288230376151711744} &\leq \langle \left(1,\,0,\,0,\,0,\,0,\,0,\,0,\,0,\,0,\,0,\,0,\,0,\,0,\,0,\,0,\,0,\,0,\,-100,\,0,\,0\right) ,x \rangle \leq \frac{582225359826457723}{576460752303423488} \\
-\frac{309485009827534768921207957}{309485009821345068724781056} &\leq \langle \left(0,\,0,\,0,\,0,\,0,\,0,\,0,\,0,\,0,\,1,\,0,\,0,\,0,\,0,\,0,\,0,\,0,\,0,\,0,\,0\right) ,x \rangle \leq \frac{309485009827534768921207957}{309485009821345068724781056} \\
-\frac{5332601763179159}{1125899906842624} &\leq \langle \left(-100,\,0,\,0,\,-16,\,15,\,0,\,0,\,0,\,0,\,0,\,0,\,0,\,0,\,0,\,0,\,0,\,-1,\,-1,\,0,\,0\right) ,x \rangle \leq 1 \\
-\frac{10016645324044103}{9007199254740992} &\leq \langle \left(0,\,0,\,0,\,1,\,0,\,0,\,0,\,0,\,0,\,0,\,0,\,0,\,0,\,0,\,0,\,0,\,0,\,-16,\,0,\,0\right) ,x \rangle \leq \frac{40081011649042577}{36028797018963968} \\
-\frac{1180593503442084915409}{1180591620717411303424} &\leq \langle \left(0,\,0,\,0,\,0,\,0,\,0,\,0,\,1,\,0,\,0,\,0,\,0,\,0,\,0,\,0,\,0,\,0,\,0,\,0,\,0\right) ,x \rangle \leq \frac{1180593503442084915409}{1180591620717411303424} \\
-\frac{151699951751958945}{144115188075855872} &\leq \langle \left(0,\,0,\,1,\,0,\,0,\,0,\,0,\,0,\,0,\,0,\,0,\,0,\,0,\,0,\,0,\,0,\,0,\,0,\,0,\,0\right) ,x \rangle \leq \frac{80781273703537389}{72057594037927936} \\
-\frac{2417851643085731713983019}{2417851639229258349412352} &\leq \langle \left(0,\,0,\,0,\,0,\,0,\,0,\,0,\,0,\,1,\,0,\,0,\,0,\,0,\,0,\,0,\,0,\,0,\,0,\,0,\,0\right) ,x \rangle \leq \frac{2417851643085731713983019}{2417851639229258349412352} \\
-\frac{2361185035008561037417}{2361183241434822606848} &\leq \langle \left(0,\,0,\,0,\,0,\,0,\,0,\,0,\,0,\,0,\,0,\,0,\,0,\,0,\,0,\,1,\,0,\,0,\,0,\,0,\,0\right) ,x \rangle \leq \frac{2361185035008561037417}{2361183241434822606848} \\
-\frac{18889473462374272866723}{18889465931478580854784} &\leq \langle \left(0,\,0,\,0,\,0,\,0,\,1,\,0,\,0,\,0,\,0,\,0,\,0,\,0,\,0,\,0,\,0,\,0,\,0,\,0,\,0\right) ,x \rangle \leq \frac{18889473462375840754951}{18889465931478580854784} \\
\end{align*}}
\subsection{Counterexample based on inequality $2$}
If $H$ is the polyhedron obtained from $N$ by removing inequality $I_{2}$, then $H$ is a bounded Hirsch polytope and $N = H \cap C$, where $C$ is the combinatorial $20$-cube defined by
{\footnotesize\begin{align*}
-\frac{1180593414291149733993}{1180591620717411303424} &\leq \langle \left(0,\,0,\,0,\,0,\,0,\,0,\,0,\,0,\,0,\,0,\,0,\,0,\,0,\,0,\,0,\,1,\,0,\,0,\,0,\,0\right) ,x \rangle \leq \frac{1180593414291149733993}{1180591620717411303424} \\
-\frac{38280228648047971}{36028797018963968} &\leq \langle \left(0,\,0,\,0,\,0,\,1,\,0,\,0,\,0,\,0,\,0,\,0,\,0,\,0,\,0,\,0,\,0,\,0,\,0,\,32,\,0\right) ,x \rangle \leq \frac{76542601517243887}{72057594037927936} \\
-\frac{2361185124158778945083}{2361183241434822606848} &\leq \langle \left(0,\,0,\,0,\,0,\,0,\,0,\,1,\,0,\,0,\,0,\,0,\,0,\,0,\,0,\,0,\,0,\,0,\,0,\,0,\,0\right) ,x \rangle \leq \frac{2361185124159496218833}{2361183241434822606848} \\
-\frac{18889470523027351237041}{18889465931478580854784} &\leq \langle \left(0,\,0,\,0,\,0,\,0,\,0,\,0,\,0,\,0,\,0,\,0,\,0,\,0,\,0,\,0,\,0,\,0,\,0,\,0,\,1\right) ,x \rangle \leq \frac{18889470523027351237041}{18889465931478580854784} \\
-\frac{76859732694334399}{72057594037927936} &\leq \langle \left(0,\,1,\,0,\,0,\,0,\,0,\,0,\,0,\,0,\,0,\,0,\,0,\,0,\,0,\,0,\,0,\,0,\,0,\,0,\,0\right) ,x \rangle \leq \frac{77296343773163779}{72057594037927936} \\
-\frac{590297603932444082281}{590295810358705651712} &\leq \langle \left(0,\,0,\,0,\,0,\,0,\,0,\,0,\,0,\,0,\,0,\,0,\,0,\,0,\,0,\,0,\,0,\,1,\,0,\,0,\,0\right) ,x \rangle \leq \frac{590297603932444082281}{590295810358705651712} \\
-\frac{39614081257140112583132697033}{39614081257132168796771975168} &\leq \langle \left(0,\,0,\,0,\,0,\,0,\,0,\,0,\,0,\,0,\,0,\,1,\,0,\,0,\,0,\,0,\,0,\,0,\,0,\,0,\,0\right) ,x \rangle \leq \frac{9903520314285028145783175131}{9903520314283042199192993792} \\
-\frac{2417851645399615732725419}{2417851639229258349412352} &\leq \langle \left(0,\,0,\,0,\,0,\,0,\,0,\,0,\,0,\,0,\,0,\,0,\,0,\,1,\,0,\,0,\,0,\,0,\,0,\,0,\,0\right) ,x \rangle \leq \frac{2417851645399615732725419}{2417851639229258349412352} \\
-\frac{147579558414782802697}{147573952589676412928} &\leq \langle \left(0,\,0,\,0,\,0,\,0,\,0,\,0,\,0,\,0,\,0,\,0,\,0,\,0,\,1,\,0,\,0,\,0,\,0,\,0,\,0\right) ,x \rangle \leq \frac{147579558414782802697}{147573952589676412928} \\
-\frac{9005724959276317}{2251799813685248} &\leq \langle \left(-100,\,0,\,0,\,0,\,32,\,0,\,0,\,0,\,0,\,0,\,0,\,0,\,0,\,0,\,0,\,0,\,0,\,-1,\,-1,\,0\right) ,x \rangle \leq 1 \\
-\frac{9671406563106733594076309}{9671406556917033397649408} &\leq \langle \left(0,\,0,\,0,\,0,\,0,\,0,\,0,\,0,\,0,\,0,\,0,\,1,\,0,\,0,\,0,\,0,\,0,\,0,\,0,\,0\right) ,x \rangle \leq \frac{9671406563106733594076309}{9671406556917033397649408} \\
-\frac{309485009827534768921207957}{309485009821345068724781056} &\leq \langle \left(0,\,0,\,0,\,0,\,0,\,0,\,0,\,0,\,0,\,1,\,0,\,0,\,0,\,0,\,0,\,0,\,0,\,0,\,0,\,0\right) ,x \rangle \leq \frac{309485009827534768921207957}{309485009821345068724781056} \\
-\frac{9444738705175253405213}{9444732965739290427392} &\leq \langle \left(0,\,0,\,0,\,0,\,0,\,0,\,0,\,0,\,0,\,0,\,0,\,0,\,0,\,0,\,0,\,0,\,0,\,1,\,-1,\,0\right) ,x \rangle \leq \frac{1180593773005897420107}{1180591620717411303424} \\
-\frac{151699951751958945}{144115188075855872} &\leq \langle \left(0,\,0,\,1,\,0,\,0,\,0,\,0,\,0,\,0,\,0,\,0,\,0,\,0,\,0,\,0,\,0,\,0,\,0,\,0,\,0\right) ,x \rangle \leq \frac{80781273703537389}{72057594037927936} \\
-\frac{1180593503442084915409}{1180591620717411303424} &\leq \langle \left(0,\,0,\,0,\,0,\,0,\,0,\,0,\,1,\,0,\,0,\,0,\,0,\,0,\,0,\,0,\,0,\,0,\,0,\,0,\,0\right) ,x \rangle \leq \frac{1180593503442084915409}{1180591620717411303424} \\
-\frac{18889473462374272866723}{18889465931478580854784} &\leq \langle \left(0,\,0,\,0,\,0,\,0,\,1,\,0,\,0,\,0,\,0,\,0,\,0,\,0,\,0,\,0,\,0,\,0,\,0,\,0,\,0\right) ,x \rangle \leq \frac{18889473462375840754951}{18889465931478580854784} \\
-\frac{2417851643085731713983019}{2417851639229258349412352} &\leq \langle \left(0,\,0,\,0,\,0,\,0,\,0,\,0,\,0,\,1,\,0,\,0,\,0,\,0,\,0,\,0,\,0,\,0,\,0,\,0,\,0\right) ,x \rangle \leq \frac{2417851643085731713983019}{2417851639229258349412352} \\
-\frac{2361185035008561037417}{2361183241434822606848} &\leq \langle \left(0,\,0,\,0,\,0,\,0,\,0,\,0,\,0,\,0,\,0,\,0,\,0,\,0,\,0,\,1,\,0,\,0,\,0,\,0,\,0\right) ,x \rangle \leq \frac{2361185035008561037417}{2361183241434822606848} \\
-\frac{291145675183409877}{288230376151711744} &\leq \langle \left(1,\,0,\,0,\,0,\,0,\,0,\,0,\,0,\,0,\,0,\,0,\,0,\,0,\,0,\,0,\,0,\,0,\,0,\,-100,\,0\right) ,x \rangle \leq \frac{582225359826457723}{576460752303423488} \\
-\frac{40066579960122721}{36028797018963968} &\leq \langle \left(0,\,0,\,0,\,1,\,0,\,0,\,0,\,0,\,0,\,0,\,0,\,0,\,0,\,0,\,0,\,0,\,0,\,0,\,0,\,0\right) ,x \rangle \leq \frac{20040505775302443}{18014398509481984} \\
\end{align*}}
\subsection{Counterexample based on inequality $3$}
If $H$ is the polyhedron obtained from $N$ by removing inequality $I_{3}$, then $H$ is a bounded Hirsch polytope and $N = H \cap C$, where $C$ is the combinatorial $20$-cube defined by
{\footnotesize\begin{align*}
-\frac{4722372222305608191517}{4722366482869645213696} &\leq \langle \left(0,\,0,\,0,\,0,\,0,\,0,\,0,\,0,\,0,\,0,\,0,\,0,\,0,\,0,\,0,\,0,\,0,\,0,\,1,\,-1\right) ,x \rangle \leq \frac{4722373370192800787081}{4722366482869645213696} \\
-\frac{9671406563106733594076309}{9671406556917033397649408} &\leq \langle \left(0,\,0,\,0,\,0,\,0,\,0,\,0,\,0,\,0,\,0,\,0,\,1,\,0,\,0,\,0,\,0,\,0,\,0,\,0,\,0\right) ,x \rangle \leq \frac{9671406563106733594076309}{9671406556917033397649408} \\
-\frac{2361185124158778945083}{2361183241434822606848} &\leq \langle \left(0,\,0,\,0,\,0,\,0,\,0,\,1,\,0,\,0,\,0,\,0,\,0,\,0,\,0,\,0,\,0,\,0,\,0,\,0,\,0\right) ,x \rangle \leq \frac{2361185124159496218833}{2361183241434822606848} \\
-\frac{76859732694334399}{72057594037927936} &\leq \langle \left(0,\,1,\,0,\,0,\,0,\,0,\,0,\,0,\,0,\,0,\,0,\,0,\,0,\,0,\,0,\,0,\,0,\,0,\,0,\,0\right) ,x \rangle \leq \frac{77296343773163779}{72057594037927936} \\
-\frac{590297603932444082281}{590295810358705651712} &\leq \langle \left(0,\,0,\,0,\,0,\,0,\,0,\,0,\,0,\,0,\,0,\,0,\,0,\,0,\,0,\,0,\,0,\,1,\,0,\,0,\,0\right) ,x \rangle \leq \frac{590297603932444082281}{590295810358705651712} \\
-\frac{9444738705175253405213}{9444732965739290427392} &\leq \langle \left(0,\,0,\,0,\,0,\,0,\,0,\,0,\,0,\,0,\,0,\,0,\,0,\,0,\,0,\,0,\,0,\,0,\,1,\,0,\,0\right) ,x \rangle \leq \frac{9444738705175253405213}{9444732965739290427392} \\
-\frac{76560456587202163}{72057594037927936} &\leq \langle \left(0,\,0,\,0,\,0,\,1,\,0,\,0,\,0,\,0,\,0,\,0,\,0,\,0,\,0,\,0,\,0,\,0,\,0,\,0,\,15\right) ,x \rangle \leq \frac{76542601517243887}{72057594037927936} \\
-\frac{39614081257140112583132697033}{39614081257132168796771975168} &\leq \langle \left(0,\,0,\,0,\,0,\,0,\,0,\,0,\,0,\,0,\,0,\,1,\,0,\,0,\,0,\,0,\,0,\,0,\,0,\,0,\,0\right) ,x \rangle \leq \frac{9903520314285028145783175131}{9903520314283042199192993792} \\
-\frac{2417851645399615732725419}{2417851639229258349412352} &\leq \langle \left(0,\,0,\,0,\,0,\,0,\,0,\,0,\,0,\,0,\,0,\,0,\,0,\,1,\,0,\,0,\,0,\,0,\,0,\,0,\,0\right) ,x \rangle \leq \frac{2417851645399615732725419}{2417851639229258349412352} \\
-\frac{147579558414782802697}{147573952589676412928} &\leq \langle \left(0,\,0,\,0,\,0,\,0,\,0,\,0,\,0,\,0,\,0,\,0,\,0,\,0,\,1,\,0,\,0,\,0,\,0,\,0,\,0\right) ,x \rangle \leq \frac{147579558414782802697}{147573952589676412928} \\
-\frac{2361185035008561037417}{2361183241434822606848} &\leq \langle \left(0,\,0,\,0,\,0,\,0,\,0,\,0,\,0,\,0,\,0,\,0,\,0,\,0,\,0,\,1,\,0,\,0,\,0,\,0,\,0\right) ,x \rangle \leq \frac{2361185035008561037417}{2361183241434822606848} \\
-\frac{582236090415071961}{576460752303423488} &\leq \langle \left(1,\,0,\,0,\,0,\,0,\,0,\,0,\,0,\,0,\,0,\,0,\,0,\,0,\,0,\,0,\,0,\,0,\,0,\,0,\,-100\right) ,x \rangle \leq \frac{582225359826457723}{576460752303423488} \\
-\frac{309485009827534768921207957}{309485009821345068724781056} &\leq \langle \left(0,\,0,\,0,\,0,\,0,\,0,\,0,\,0,\,0,\,1,\,0,\,0,\,0,\,0,\,0,\,0,\,0,\,0,\,0,\,0\right) ,x \rangle \leq \frac{309485009827534768921207957}{309485009821345068724781056} \\
-\frac{10627220392380869}{2251799813685248} &\leq \langle \left(-100,\,0,\,0,\,16,\,15,\,0,\,0,\,0,\,0,\,0,\,0,\,0,\,0,\,0,\,0,\,0,\,0,\,0,\,-1,\,-1\right) ,x \rangle \leq 1 \\
-\frac{1180593414291149733993}{1180591620717411303424} &\leq \langle \left(0,\,0,\,0,\,0,\,0,\,0,\,0,\,0,\,0,\,0,\,0,\,0,\,0,\,0,\,0,\,1,\,0,\,0,\,0,\,0\right) ,x \rangle \leq \frac{1180593414291149733993}{1180591620717411303424} \\
-\frac{1180593503442084915409}{1180591620717411303424} &\leq \langle \left(0,\,0,\,0,\,0,\,0,\,0,\,0,\,1,\,0,\,0,\,0,\,0,\,0,\,0,\,0,\,0,\,0,\,0,\,0,\,0\right) ,x \rangle \leq \frac{1180593503442084915409}{1180591620717411303424} \\
-\frac{151699951751958945}{144115188075855872} &\leq \langle \left(0,\,0,\,1,\,0,\,0,\,0,\,0,\,0,\,0,\,0,\,0,\,0,\,0,\,0,\,0,\,0,\,0,\,0,\,0,\,0\right) ,x \rangle \leq \frac{80781273703537389}{72057594037927936} \\
-\frac{2417851643085731713983019}{2417851639229258349412352} &\leq \langle \left(0,\,0,\,0,\,0,\,0,\,0,\,0,\,0,\,1,\,0,\,0,\,0,\,0,\,0,\,0,\,0,\,0,\,0,\,0,\,0\right) ,x \rangle \leq \frac{2417851643085731713983019}{2417851639229258349412352} \\
-\frac{20033290247272103}{18014398509481984} &\leq \langle \left(0,\,0,\,0,\,1,\,0,\,0,\,0,\,0,\,0,\,0,\,0,\,0,\,0,\,0,\,0,\,0,\,0,\,0,\,0,\,16\right) ,x \rangle \leq \frac{10020252897507271}{9007199254740992} \\
-\frac{18889473462374272866723}{18889465931478580854784} &\leq \langle \left(0,\,0,\,0,\,0,\,0,\,1,\,0,\,0,\,0,\,0,\,0,\,0,\,0,\,0,\,0,\,0,\,0,\,0,\,0,\,0\right) ,x \rangle \leq \frac{18889473462375840754951}{18889465931478580854784} \\
\end{align*}}
\subsection{Counterexample based on inequality $4$}
If $H$ is the polyhedron obtained from $N$ by removing inequality $I_{4}$, then $H$ is a bounded Hirsch polytope and $N = H \cap C$, where $C$ is the combinatorial $20$-cube defined by
{\footnotesize\begin{align*}
-\frac{4041417668546867}{562949953421312} &\leq \langle \left(7,\,0,\,0,\,0,\,100,\,0,\,0,\,0,\,0,\,0,\,0,\,0,\,0,\,0,\,0,\,0,\,0,\,0,\,0,\,0\right) ,x \rangle \leq \frac{1026485795196359}{140737488355328} \\
-\frac{1180593503442084915409}{1180591620717411303424} &\leq \langle \left(0,\,0,\,0,\,0,\,0,\,0,\,0,\,1,\,0,\,0,\,0,\,0,\,0,\,0,\,0,\,0,\,0,\,0,\,0,\,0\right) ,x \rangle \leq \frac{1180593503442084915409}{1180591620717411303424} \\
-\frac{4722372222305608191517}{4722366482869645213696} &\leq \langle \left(0,\,0,\,0,\,0,\,0,\,0,\,0,\,0,\,0,\,0,\,0,\,0,\,0,\,0,\,0,\,0,\,0,\,0,\,1,\,0\right) ,x \rangle \leq \frac{4722372222305608191517}{4722366482869645213696} \\
-\frac{76859732694334399}{72057594037927936} &\leq \langle \left(0,\,1,\,0,\,0,\,0,\,0,\,0,\,0,\,0,\,0,\,0,\,0,\,0,\,0,\,0,\,0,\,0,\,0,\,0,\,0\right) ,x \rangle \leq \frac{77296343773163779}{72057594037927936} \\
-\frac{18889470523027351237041}{18889465931478580854784} &\leq \langle \left(0,\,0,\,0,\,0,\,0,\,0,\,0,\,0,\,0,\,0,\,0,\,0,\,0,\,0,\,0,\,0,\,0,\,0,\,0,\,1\right) ,x \rangle \leq \frac{18889470523027351237041}{18889465931478580854784} \\
-\frac{2417851645399615732725419}{2417851639229258349412352} &\leq \langle \left(0,\,0,\,0,\,0,\,0,\,0,\,0,\,0,\,0,\,0,\,0,\,0,\,1,\,0,\,0,\,0,\,0,\,0,\,0,\,0\right) ,x \rangle \leq \frac{2417851645399615732725419}{2417851639229258349412352} \\
-\frac{590297603932444082281}{590295810358705651712} &\leq \langle \left(0,\,0,\,0,\,0,\,0,\,0,\,0,\,0,\,0,\,0,\,0,\,0,\,0,\,0,\,0,\,0,\,1,\,0,\,0,\,0\right) ,x \rangle \leq \frac{590297603932444082281}{590295810358705651712} \\
-\frac{9444738705175253405213}{9444732965739290427392} &\leq \langle \left(0,\,0,\,0,\,0,\,0,\,0,\,0,\,0,\,0,\,0,\,0,\,0,\,0,\,0,\,0,\,0,\,0,\,1,\,0,\,0\right) ,x \rangle \leq \frac{9444738705175253405213}{9444732965739290427392} \\
-\frac{39614081257140112583132697033}{39614081257132168796771975168} &\leq \langle \left(0,\,0,\,0,\,0,\,0,\,0,\,0,\,0,\,0,\,0,\,1,\,0,\,0,\,0,\,0,\,0,\,0,\,0,\,0,\,0\right) ,x \rangle \leq \frac{9903520314285028145783175131}{9903520314283042199192993792} \\
-\frac{11271016221296233}{9007199254740992} &\leq \langle \left(0,\,0,\,0,\,1,\,-3,\,0,\,0,\,0,\,0,\,0,\,0,\,0,\,0,\,0,\,0,\,0,\,0,\,0,\,0,\,0\right) ,x \rangle \leq \frac{11708675609116937}{9007199254740992} \\
-\frac{147579558414782802697}{147573952589676412928} &\leq \langle \left(0,\,0,\,0,\,0,\,0,\,0,\,0,\,0,\,0,\,0,\,0,\,0,\,0,\,1,\,0,\,0,\,0,\,0,\,0,\,0\right) ,x \rangle \leq \frac{147579558414782802697}{147573952589676412928} \\
-\frac{9671406563106733594076309}{9671406556917033397649408} &\leq \langle \left(0,\,0,\,0,\,0,\,0,\,0,\,0,\,0,\,0,\,0,\,0,\,1,\,0,\,0,\,0,\,0,\,0,\,0,\,0,\,0\right) ,x \rangle \leq \frac{9671406563106733594076309}{9671406556917033397649408} \\
-\frac{18889473462374272866723}{18889465931478580854784} &\leq \langle \left(0,\,0,\,0,\,0,\,0,\,1,\,0,\,0,\,0,\,0,\,0,\,0,\,0,\,0,\,0,\,0,\,0,\,0,\,0,\,0\right) ,x \rangle \leq \frac{18889473462375840754951}{18889465931478580854784} \\
-\frac{309485009827534768921207957}{309485009821345068724781056} &\leq \langle \left(0,\,0,\,0,\,0,\,0,\,0,\,0,\,0,\,0,\,1,\,0,\,0,\,0,\,0,\,0,\,0,\,0,\,0,\,0,\,0\right) ,x \rangle \leq \frac{309485009827534768921207957}{309485009821345068724781056} \\
-\frac{1180593414291149733993}{1180591620717411303424} &\leq \langle \left(0,\,0,\,0,\,0,\,0,\,0,\,0,\,0,\,0,\,0,\,0,\,0,\,0,\,0,\,0,\,1,\,0,\,0,\,0,\,0\right) ,x \rangle \leq \frac{1180593414291149733993}{1180591620717411303424} \\
-\frac{2361185124158778945083}{2361183241434822606848} &\leq \langle \left(0,\,0,\,0,\,0,\,0,\,0,\,1,\,0,\,0,\,0,\,0,\,0,\,0,\,0,\,0,\,0,\,0,\,0,\,0,\,0\right) ,x \rangle \leq \frac{2361185124159496218833}{2361183241434822606848} \\
-\frac{2417851643085731713983019}{2417851639229258349412352} &\leq \langle \left(0,\,0,\,0,\,0,\,0,\,0,\,0,\,0,\,1,\,0,\,0,\,0,\,0,\,0,\,0,\,0,\,0,\,0,\,0,\,0\right) ,x \rangle \leq \frac{2417851643085731713983019}{2417851639229258349412352} \\
-\frac{2361185035008561037417}{2361183241434822606848} &\leq \langle \left(0,\,0,\,0,\,0,\,0,\,0,\,0,\,0,\,0,\,0,\,0,\,0,\,0,\,0,\,1,\,0,\,0,\,0,\,0,\,0\right) ,x \rangle \leq \frac{2361185035008561037417}{2361183241434822606848} \\
-\frac{10767421176423341}{2251799813685248} &\leq \langle \left(-100,\,0,\,0,\,21,\,7,\,0,\,0,\,0,\,0,\,0,\,0,\,0,\,0,\,0,\,0,\,0,\,0,\,0,\,0,\,0\right) ,x \rangle \leq 1 \\
-\frac{151699951751958945}{144115188075855872} &\leq \langle \left(0,\,0,\,1,\,0,\,0,\,0,\,0,\,0,\,0,\,0,\,0,\,0,\,0,\,0,\,0,\,0,\,0,\,0,\,0,\,0\right) ,x \rangle \leq \frac{80781273703537389}{72057594037927936} \\
\end{align*}}
\subsection{Counterexample based on inequality $5$}
If $H$ is the polyhedron obtained from $N$ by removing inequality $I_{5}$, then $H$ is a bounded Hirsch polytope and $N = H \cap C$, where $C$ is the combinatorial $20$-cube defined by
{\footnotesize\begin{align*}
-\frac{5047228677492251}{1125899906842624} &\leq \langle \left(-100,\,0,\,0,\,20,\,4,\,0,\,0,\,0,\,0,\,0,\,0,\,0,\,0,\,0,\,0,\,0,\,0,\,0,\,0,\,-1\right) ,x \rangle \leq 1 \\
-\frac{147579558414782802697}{147573952589676412928} &\leq \langle \left(0,\,0,\,0,\,0,\,0,\,0,\,0,\,0,\,0,\,0,\,0,\,0,\,0,\,1,\,0,\,0,\,0,\,0,\,0,\,0\right) ,x \rangle \leq \frac{147579558414782802697}{147573952589676412928} \\
-\frac{80133161256299147}{72057594037927936} &\leq \langle \left(0,\,0,\,0,\,1,\,0,\,0,\,0,\,0,\,0,\,0,\,0,\,0,\,0,\,0,\,0,\,0,\,0,\,0,\,0,\,20\right) ,x \rangle \leq \frac{80162023199770267}{72057594037927936} \\
-\frac{2361185124158778945083}{2361183241434822606848} &\leq \langle \left(0,\,0,\,0,\,0,\,0,\,0,\,1,\,0,\,0,\,0,\,0,\,0,\,0,\,0,\,0,\,0,\,0,\,0,\,0,\,0\right) ,x \rangle \leq \frac{2361185124159496218833}{2361183241434822606848} \\
-\frac{4722372222305608191517}{4722366482869645213696} &\leq \langle \left(0,\,0,\,0,\,0,\,0,\,0,\,0,\,0,\,0,\,0,\,0,\,0,\,0,\,0,\,0,\,0,\,0,\,0,\,1,\,0\right) ,x \rangle \leq \frac{4722372222305608191517}{4722366482869645213696} \\
-\frac{76859732694334399}{72057594037927936} &\leq \langle \left(0,\,1,\,0,\,0,\,0,\,0,\,0,\,0,\,0,\,0,\,0,\,0,\,0,\,0,\,0,\,0,\,0,\,0,\,0,\,0\right) ,x \rangle \leq \frac{77296343773163779}{72057594037927936} \\
-\frac{590297603932444082281}{590295810358705651712} &\leq \langle \left(0,\,0,\,0,\,0,\,0,\,0,\,0,\,0,\,0,\,0,\,0,\,0,\,0,\,0,\,0,\,0,\,1,\,0,\,0,\,0\right) ,x \rangle \leq \frac{590297603932444082281}{590295810358705651712} \\
-\frac{9444738705175253405213}{9444732965739290427392} &\leq \langle \left(0,\,0,\,0,\,0,\,0,\,0,\,0,\,0,\,0,\,0,\,0,\,0,\,0,\,0,\,0,\,0,\,0,\,1,\,0,\,0\right) ,x \rangle \leq \frac{9444738705175253405213}{9444732965739290427392} \\
-\frac{39614081257140112583132697033}{39614081257132168796771975168} &\leq \langle \left(0,\,0,\,0,\,0,\,0,\,0,\,0,\,0,\,0,\,0,\,1,\,0,\,0,\,0,\,0,\,0,\,0,\,0,\,0,\,0\right) ,x \rangle \leq \frac{9903520314285028145783175131}{9903520314283042199192993792} \\
-\frac{2417851645399615732725419}{2417851639229258349412352} &\leq \langle \left(0,\,0,\,0,\,0,\,0,\,0,\,0,\,0,\,0,\,0,\,0,\,0,\,1,\,0,\,0,\,0,\,0,\,0,\,0,\,0\right) ,x \rangle \leq \frac{2417851645399615732725419}{2417851639229258349412352} \\
-\frac{153120913066848029}{144115188075855872} &\leq \langle \left(0,\,0,\,0,\,0,\,1,\,0,\,0,\,0,\,0,\,0,\,0,\,0,\,0,\,0,\,0,\,0,\,0,\,0,\,0,\,4\right) ,x \rangle \leq \frac{76542601517243887}{72057594037927936} \\
-\frac{9671406563106733594076309}{9671406556917033397649408} &\leq \langle \left(0,\,0,\,0,\,0,\,0,\,0,\,0,\,0,\,0,\,0,\,0,\,1,\,0,\,0,\,0,\,0,\,0,\,0,\,0,\,0\right) ,x \rangle \leq \frac{9671406563106733594076309}{9671406556917033397649408} \\
-\frac{582236090415071961}{576460752303423488} &\leq \langle \left(1,\,0,\,0,\,0,\,0,\,0,\,0,\,0,\,0,\,0,\,0,\,0,\,0,\,0,\,0,\,0,\,0,\,0,\,0,\,-100\right) ,x \rangle \leq \frac{582225359826457723}{576460752303423488} \\
-\frac{309485009827534768921207957}{309485009821345068724781056} &\leq \langle \left(0,\,0,\,0,\,0,\,0,\,0,\,0,\,0,\,0,\,1,\,0,\,0,\,0,\,0,\,0,\,0,\,0,\,0,\,0,\,0\right) ,x \rangle \leq \frac{309485009827534768921207957}{309485009821345068724781056} \\
-\frac{1180593414291149733993}{1180591620717411303424} &\leq \langle \left(0,\,0,\,0,\,0,\,0,\,0,\,0,\,0,\,0,\,0,\,0,\,0,\,0,\,0,\,0,\,1,\,0,\,0,\,0,\,0\right) ,x \rangle \leq \frac{1180593414291149733993}{1180591620717411303424} \\
-\frac{1180593503442084915409}{1180591620717411303424} &\leq \langle \left(0,\,0,\,0,\,0,\,0,\,0,\,0,\,1,\,0,\,0,\,0,\,0,\,0,\,0,\,0,\,0,\,0,\,0,\,0,\,0\right) ,x \rangle \leq \frac{1180593503442084915409}{1180591620717411303424} \\
-\frac{151699951751958945}{144115188075855872} &\leq \langle \left(0,\,0,\,1,\,0,\,0,\,0,\,0,\,0,\,0,\,0,\,0,\,0,\,0,\,0,\,0,\,0,\,0,\,0,\,0,\,0\right) ,x \rangle \leq \frac{80781273703537389}{72057594037927936} \\
-\frac{2417851643085731713983019}{2417851639229258349412352} &\leq \langle \left(0,\,0,\,0,\,0,\,0,\,0,\,0,\,0,\,1,\,0,\,0,\,0,\,0,\,0,\,0,\,0,\,0,\,0,\,0,\,0\right) ,x \rangle \leq \frac{2417851643085731713983019}{2417851639229258349412352} \\
-\frac{2361185035008561037417}{2361183241434822606848} &\leq \langle \left(0,\,0,\,0,\,0,\,0,\,0,\,0,\,0,\,0,\,0,\,0,\,0,\,0,\,0,\,1,\,0,\,0,\,0,\,0,\,0\right) ,x \rangle \leq \frac{2361185035008561037417}{2361183241434822606848} \\
-\frac{18889473462374272866723}{18889465931478580854784} &\leq \langle \left(0,\,0,\,0,\,0,\,0,\,1,\,0,\,0,\,0,\,0,\,0,\,0,\,0,\,0,\,0,\,0,\,0,\,0,\,0,\,0\right) ,x \rangle \leq \frac{18889473462375840754951}{18889465931478580854784} \\
\end{align*}}
\subsection{Counterexample based on inequality $6$}
If $H$ is the polyhedron obtained from $N$ by removing inequality $I_{6}$, then $H$ is a bounded Hirsch polytope and $N = H \cap C$, where $C$ is the combinatorial $20$-cube defined by
{\footnotesize\begin{align*}
-\frac{4485007479315951}{1125899906842624} &\leq \langle \left(-100,\,0,\,0,\,0,\,-32,\,0,\,0,\,0,\,0,\,0,\,0,\,0,\,0,\,0,\,0,\,0,\,0,\,0,\,0,\,0\right) ,x \rangle \leq 1 \\
-\frac{9671406563106733594076309}{9671406556917033397649408} &\leq \langle \left(0,\,0,\,0,\,0,\,0,\,0,\,0,\,0,\,0,\,0,\,0,\,1,\,0,\,0,\,0,\,0,\,0,\,0,\,0,\,0\right) ,x \rangle \leq \frac{9671406563106733594076309}{9671406556917033397649408} \\
-\frac{2361185124158778945083}{2361183241434822606848} &\leq \langle \left(0,\,0,\,0,\,0,\,0,\,0,\,1,\,0,\,0,\,0,\,0,\,0,\,0,\,0,\,0,\,0,\,0,\,0,\,0,\,0\right) ,x \rangle \leq \frac{2361185124159496218833}{2361183241434822606848} \\
-\frac{4722372222305608191517}{4722366482869645213696} &\leq \langle \left(0,\,0,\,0,\,0,\,0,\,0,\,0,\,0,\,0,\,0,\,0,\,0,\,0,\,0,\,0,\,0,\,0,\,0,\,1,\,0\right) ,x \rangle \leq \frac{4722372222305608191517}{4722366482869645213696} \\
-\frac{76859732694334399}{72057594037927936} &\leq \langle \left(0,\,1,\,0,\,0,\,0,\,0,\,0,\,0,\,0,\,0,\,0,\,0,\,0,\,0,\,0,\,0,\,0,\,0,\,0,\,0\right) ,x \rangle \leq \frac{77296343773163779}{72057594037927936} \\
-\frac{18889470523027351237041}{18889465931478580854784} &\leq \langle \left(0,\,0,\,0,\,0,\,0,\,0,\,0,\,0,\,0,\,0,\,0,\,0,\,0,\,0,\,0,\,0,\,0,\,0,\,0,\,1\right) ,x \rangle \leq \frac{18889470523027351237041}{18889465931478580854784} \\
-\frac{590297603932444082281}{590295810358705651712} &\leq \langle \left(0,\,0,\,0,\,0,\,0,\,0,\,0,\,0,\,0,\,0,\,0,\,0,\,0,\,0,\,0,\,0,\,1,\,0,\,0,\,0\right) ,x \rangle \leq \frac{590297603932444082281}{590295810358705651712} \\
-\frac{9444738705175253405213}{9444732965739290427392} &\leq \langle \left(0,\,0,\,0,\,0,\,0,\,0,\,0,\,0,\,0,\,0,\,0,\,0,\,0,\,0,\,0,\,0,\,0,\,1,\,0,\,0\right) ,x \rangle \leq \frac{9444738705175253405213}{9444732965739290427392} \\
-\frac{39614081257140112583132697033}{39614081257132168796771975168} &\leq \langle \left(0,\,0,\,0,\,0,\,0,\,0,\,0,\,0,\,0,\,0,\,1,\,0,\,0,\,0,\,0,\,0,\,0,\,0,\,0,\,0\right) ,x \rangle \leq \frac{9903520314285028145783175131}{9903520314283042199192993792} \\
-\frac{2417851645399615732725419}{2417851639229258349412352} &\leq \langle \left(0,\,0,\,0,\,0,\,0,\,0,\,0,\,0,\,0,\,0,\,0,\,0,\,1,\,0,\,0,\,0,\,0,\,0,\,0,\,0\right) ,x \rangle \leq \frac{2417851645399615732725419}{2417851639229258349412352} \\
-\frac{147579558414782802697}{147573952589676412928} &\leq \langle \left(0,\,0,\,0,\,0,\,0,\,0,\,0,\,0,\,0,\,0,\,0,\,0,\,0,\,1,\,0,\,0,\,0,\,0,\,0,\,0\right) ,x \rangle \leq \frac{147579558414782802697}{147573952589676412928} \\
-\frac{11154110587300757}{4503599627370496} &\leq \langle \left(8,\,0,\,0,\,0,\,-25,\,0,\,0,\,0,\,0,\,0,\,0,\,0,\,0,\,0,\,0,\,0,\,0,\,0,\,0,\,0\right) ,x \rangle \leq \frac{11899492277370701}{4503599627370496} \\
-\frac{18889473462374272866723}{18889465931478580854784} &\leq \langle \left(0,\,0,\,0,\,0,\,0,\,1,\,0,\,0,\,0,\,0,\,0,\,0,\,0,\,0,\,0,\,0,\,0,\,0,\,0,\,0\right) ,x \rangle \leq \frac{18889473462375840754951}{18889465931478580854784} \\
-\frac{309485009827534768921207957}{309485009821345068724781056} &\leq \langle \left(0,\,0,\,0,\,0,\,0,\,0,\,0,\,0,\,0,\,1,\,0,\,0,\,0,\,0,\,0,\,0,\,0,\,0,\,0,\,0\right) ,x \rangle \leq \frac{309485009827534768921207957}{309485009821345068724781056} \\
-\frac{1180593414291149733993}{1180591620717411303424} &\leq \langle \left(0,\,0,\,0,\,0,\,0,\,0,\,0,\,0,\,0,\,0,\,0,\,0,\,0,\,0,\,0,\,1,\,0,\,0,\,0,\,0\right) ,x \rangle \leq \frac{1180593414291149733993}{1180591620717411303424} \\
-\frac{1180593503442084915409}{1180591620717411303424} &\leq \langle \left(0,\,0,\,0,\,0,\,0,\,0,\,0,\,1,\,0,\,0,\,0,\,0,\,0,\,0,\,0,\,0,\,0,\,0,\,0,\,0\right) ,x \rangle \leq \frac{1180593503442084915409}{1180591620717411303424} \\
-\frac{151699951751958945}{144115188075855872} &\leq \langle \left(0,\,0,\,1,\,0,\,0,\,0,\,0,\,0,\,0,\,0,\,0,\,0,\,0,\,0,\,0,\,0,\,0,\,0,\,0,\,0\right) ,x \rangle \leq \frac{80781273703537389}{72057594037927936} \\
-\frac{2417851643085731713983019}{2417851639229258349412352} &\leq \langle \left(0,\,0,\,0,\,0,\,0,\,0,\,0,\,0,\,1,\,0,\,0,\,0,\,0,\,0,\,0,\,0,\,0,\,0,\,0,\,0\right) ,x \rangle \leq \frac{2417851643085731713983019}{2417851639229258349412352} \\
-\frac{2361185035008561037417}{2361183241434822606848} &\leq \langle \left(0,\,0,\,0,\,0,\,0,\,0,\,0,\,0,\,0,\,0,\,0,\,0,\,0,\,0,\,1,\,0,\,0,\,0,\,0,\,0\right) ,x \rangle \leq \frac{2361185035008561037417}{2361183241434822606848} \\
-\frac{40066579960122721}{36028797018963968} &\leq \langle \left(0,\,0,\,0,\,1,\,0,\,0,\,0,\,0,\,0,\,0,\,0,\,0,\,0,\,0,\,0,\,0,\,0,\,0,\,0,\,0\right) ,x \rangle \leq \frac{20040505775302443}{18014398509481984} \\
\end{align*}}
\subsection{Counterexample based on inequality $7$}
If $H$ is the polyhedron obtained from $N$ by removing inequality $I_{7}$, then $H$ is a bounded Hirsch polytope and $N = H \cap C$, where $C$ is the combinatorial $20$-cube defined by
{\footnotesize\begin{align*}
-\frac{147579558414782802697}{147573952589676412928} &\leq \langle \left(0,\,0,\,0,\,0,\,0,\,0,\,0,\,0,\,0,\,0,\,0,\,0,\,0,\,1,\,0,\,0,\,0,\,0,\,0,\,0\right) ,x \rangle \leq \frac{147579558414782802697}{147573952589676412928} \\
-\frac{9671406563106733594076309}{9671406556917033397649408} &\leq \langle \left(0,\,0,\,0,\,0,\,0,\,0,\,0,\,0,\,0,\,0,\,0,\,1,\,0,\,0,\,0,\,0,\,0,\,0,\,0,\,0\right) ,x \rangle \leq \frac{9671406563106733594076309}{9671406556917033397649408} \\
-\frac{1180593503442084915409}{1180591620717411303424} &\leq \langle \left(0,\,0,\,0,\,0,\,0,\,0,\,0,\,1,\,0,\,0,\,0,\,0,\,0,\,0,\,0,\,0,\,0,\,0,\,0,\,0\right) ,x \rangle \leq \frac{1180593503442084915409}{1180591620717411303424} \\
-\frac{4722372222305608191517}{4722366482869645213696} &\leq \langle \left(0,\,0,\,0,\,0,\,0,\,0,\,0,\,0,\,0,\,0,\,0,\,0,\,0,\,0,\,0,\,0,\,0,\,0,\,1,\,0\right) ,x \rangle \leq \frac{4722372222305608191517}{4722366482869645213696} \\
-\frac{76859732694334399}{72057594037927936} &\leq \langle \left(0,\,1,\,0,\,0,\,0,\,0,\,0,\,0,\,0,\,0,\,0,\,0,\,0,\,0,\,0,\,0,\,0,\,0,\,0,\,0\right) ,x \rangle \leq \frac{77296343773163779}{72057594037927936} \\
-\frac{18889470523027351237041}{18889465931478580854784} &\leq \langle \left(0,\,0,\,0,\,0,\,0,\,0,\,0,\,0,\,0,\,0,\,0,\,0,\,0,\,0,\,0,\,0,\,0,\,0,\,0,\,1\right) ,x \rangle \leq \frac{18889470523027351237041}{18889465931478580854784} \\
-\frac{590297603932444082281}{590295810358705651712} &\leq \langle \left(0,\,0,\,0,\,0,\,0,\,0,\,0,\,0,\,0,\,0,\,0,\,0,\,0,\,0,\,0,\,0,\,1,\,0,\,0,\,0\right) ,x \rangle \leq \frac{590297603932444082281}{590295810358705651712} \\
-\frac{9444738705175253405213}{9444732965739290427392} &\leq \langle \left(0,\,0,\,0,\,0,\,0,\,0,\,0,\,0,\,0,\,0,\,0,\,0,\,0,\,0,\,0,\,0,\,0,\,1,\,0,\,0\right) ,x \rangle \leq \frac{9444738705175253405213}{9444732965739290427392} \\
-\frac{39614081257140112583132697033}{39614081257132168796771975168} &\leq \langle \left(0,\,0,\,0,\,0,\,0,\,0,\,0,\,0,\,0,\,0,\,1,\,0,\,0,\,0,\,0,\,0,\,0,\,0,\,0,\,0\right) ,x \rangle \leq \frac{9903520314285028145783175131}{9903520314283042199192993792} \\
-\frac{2417851645399615732725419}{2417851639229258349412352} &\leq \langle \left(0,\,0,\,0,\,0,\,0,\,0,\,0,\,0,\,0,\,0,\,0,\,0,\,1,\,0,\,0,\,0,\,0,\,0,\,0,\,0\right) ,x \rangle \leq \frac{2417851645399615732725419}{2417851639229258349412352} \\
-\frac{5747150498801153}{2251799813685248} &\leq \langle \left(1,\,0,\,0,\,0,\,25,\,0,\,0,\,0,\,0,\,0,\,0,\,0,\,0,\,0,\,0,\,0,\,0,\,0,\,0,\,0\right) ,x \rangle \leq \frac{11556087967114283}{4503599627370496} \\
-\frac{12820660921390925}{9007199254740992} &\leq \langle \left(0,\,0,\,0,\,1,\,5,\,0,\,0,\,0,\,0,\,0,\,0,\,0,\,0,\,0,\,0,\,0,\,0,\,0,\,0,\,0\right) ,x \rangle \leq \frac{12403038523674617}{9007199254740992} \\
-\frac{18889473462374272866723}{18889465931478580854784} &\leq \langle \left(0,\,0,\,0,\,0,\,0,\,1,\,0,\,0,\,0,\,0,\,0,\,0,\,0,\,0,\,0,\,0,\,0,\,0,\,0,\,0\right) ,x \rangle \leq \frac{18889473462375840754951}{18889465931478580854784} \\
-\frac{309485009827534768921207957}{309485009821345068724781056} &\leq \langle \left(0,\,0,\,0,\,0,\,0,\,0,\,0,\,0,\,0,\,1,\,0,\,0,\,0,\,0,\,0,\,0,\,0,\,0,\,0,\,0\right) ,x \rangle \leq \frac{309485009827534768921207957}{309485009821345068724781056} \\
-\frac{1180593414291149733993}{1180591620717411303424} &\leq \langle \left(0,\,0,\,0,\,0,\,0,\,0,\,0,\,0,\,0,\,0,\,0,\,0,\,0,\,0,\,0,\,1,\,0,\,0,\,0,\,0\right) ,x \rangle \leq \frac{1180593414291149733993}{1180591620717411303424} \\
-\frac{1266317041118537}{281474976710656} &\leq \langle \left(-100,\,0,\,0,\,-20,\,4,\,0,\,0,\,0,\,0,\,0,\,0,\,0,\,0,\,0,\,0,\,0,\,0,\,0,\,0,\,0\right) ,x \rangle \leq 1 \\
-\frac{2361185124158778945083}{2361183241434822606848} &\leq \langle \left(0,\,0,\,0,\,0,\,0,\,0,\,1,\,0,\,0,\,0,\,0,\,0,\,0,\,0,\,0,\,0,\,0,\,0,\,0,\,0\right) ,x \rangle \leq \frac{2361185124159496218833}{2361183241434822606848} \\
-\frac{2417851643085731713983019}{2417851639229258349412352} &\leq \langle \left(0,\,0,\,0,\,0,\,0,\,0,\,0,\,0,\,1,\,0,\,0,\,0,\,0,\,0,\,0,\,0,\,0,\,0,\,0,\,0\right) ,x \rangle \leq \frac{2417851643085731713983019}{2417851639229258349412352} \\
-\frac{2361185035008561037417}{2361183241434822606848} &\leq \langle \left(0,\,0,\,0,\,0,\,0,\,0,\,0,\,0,\,0,\,0,\,0,\,0,\,0,\,0,\,1,\,0,\,0,\,0,\,0,\,0\right) ,x \rangle \leq \frac{2361185035008561037417}{2361183241434822606848} \\
-\frac{151699951751958945}{144115188075855872} &\leq \langle \left(0,\,0,\,1,\,0,\,0,\,0,\,0,\,0,\,0,\,0,\,0,\,0,\,0,\,0,\,0,\,0,\,0,\,0,\,0,\,0\right) ,x \rangle \leq \frac{80781273703537389}{72057594037927936} \\
\end{align*}}
\subsection{Counterexample based on inequality $8$}
If $H$ is the polyhedron obtained from $N$ by removing inequality $I_{8}$, then $H$ is a bounded Hirsch polytope and $N = H \cap C$, where $C$ is the combinatorial $20$-cube defined by
{\footnotesize\begin{align*}
-\frac{18962493993892025}{18014398509481984} &\leq \langle \left(0,\,0,\,1,\,0,\,-\frac{1}{1500},\,0,\,0,\,0,\,0,\,0,\,0,\,0,\,0,\,0,\,0,\,0,\,0,\,0,\,0,\,0\right) ,x \rangle \leq \frac{80781275884457665}{72057594037927936} \\
-\frac{1180593503442084915409}{1180591620717411303424} &\leq \langle \left(0,\,0,\,0,\,0,\,0,\,0,\,0,\,1,\,0,\,0,\,0,\,0,\,0,\,0,\,0,\,0,\,0,\,0,\,0,\,0\right) ,x \rangle \leq \frac{1180593503442084915409}{1180591620717411303424} \\
-\frac{4722372222305608191517}{4722366482869645213696} &\leq \langle \left(0,\,0,\,0,\,0,\,0,\,0,\,0,\,0,\,0,\,0,\,0,\,0,\,0,\,0,\,0,\,0,\,0,\,0,\,1,\,0\right) ,x \rangle \leq \frac{4722372222305608191517}{4722366482869645213696} \\
-\frac{4362148554682285}{1125899906842624} &\leq \langle \left(-100,\,-\frac{3}{100},\,\frac{1}{50},\,0,\,30,\,0,\,0,\,0,\,0,\,0,\,0,\,0,\,0,\,0,\,0,\,0,\,0,\,0,\,0,\,0\right) ,x \rangle \leq 1 \\
-\frac{18889470523027351237041}{18889465931478580854784} &\leq \langle \left(0,\,0,\,0,\,0,\,0,\,0,\,0,\,0,\,0,\,0,\,0,\,0,\,0,\,0,\,0,\,0,\,0,\,0,\,0,\,1\right) ,x \rangle \leq \frac{18889470523027351237041}{18889465931478580854784} \\
-\frac{590297603932444082281}{590295810358705651712} &\leq \langle \left(0,\,0,\,0,\,0,\,0,\,0,\,0,\,0,\,0,\,0,\,0,\,0,\,0,\,0,\,0,\,0,\,1,\,0,\,0,\,0\right) ,x \rangle \leq \frac{590297603932444082281}{590295810358705651712} \\
-\frac{9444738705175253405213}{9444732965739290427392} &\leq \langle \left(0,\,0,\,0,\,0,\,0,\,0,\,0,\,0,\,0,\,0,\,0,\,0,\,0,\,0,\,0,\,0,\,0,\,1,\,0,\,0\right) ,x \rangle \leq \frac{9444738705175253405213}{9444732965739290427392} \\
-\frac{39614081257140112583132697033}{39614081257132168796771975168} &\leq \langle \left(0,\,0,\,0,\,0,\,0,\,0,\,0,\,0,\,0,\,0,\,1,\,0,\,0,\,0,\,0,\,0,\,0,\,0,\,0,\,0\right) ,x \rangle \leq \frac{9903520314285028145783175131}{9903520314283042199192993792} \\
-\frac{2417851645399615732725419}{2417851639229258349412352} &\leq \langle \left(0,\,0,\,0,\,0,\,0,\,0,\,0,\,0,\,0,\,0,\,0,\,0,\,1,\,0,\,0,\,0,\,0,\,0,\,0,\,0\right) ,x \rangle \leq \frac{2417851645399615732725419}{2417851639229258349412352} \\
-\frac{147579558414782802697}{147573952589676412928} &\leq \langle \left(0,\,0,\,0,\,0,\,0,\,0,\,0,\,0,\,0,\,0,\,0,\,0,\,0,\,1,\,0,\,0,\,0,\,0,\,0,\,0\right) ,x \rangle \leq \frac{147579558414782802697}{147573952589676412928} \\
-\frac{1180593414291149733993}{1180591620717411303424} &\leq \langle \left(0,\,0,\,0,\,0,\,0,\,0,\,0,\,0,\,0,\,0,\,0,\,0,\,0,\,0,\,0,\,1,\,0,\,0,\,0,\,0\right) ,x \rangle \leq \frac{1180593414291149733993}{1180591620717411303424} \\
-\frac{9671406563106733594076309}{9671406556917033397649408} &\leq \langle \left(0,\,0,\,0,\,0,\,0,\,0,\,0,\,0,\,0,\,0,\,0,\,1,\,0,\,0,\,0,\,0,\,0,\,0,\,0,\,0\right) ,x \rangle \leq \frac{9671406563106733594076309}{9671406556917033397649408} \\
-\frac{309485009827534768921207957}{309485009821345068724781056} &\leq \langle \left(0,\,0,\,0,\,0,\,0,\,0,\,0,\,0,\,0,\,1,\,0,\,0,\,0,\,0,\,0,\,0,\,0,\,0,\,0,\,0\right) ,x \rangle \leq \frac{309485009827534768921207957}{309485009821345068724781056} \\
-\frac{674584340409955}{562949953421312} &\leq \langle \left(1,\,0,\,0,\,0,\,\frac{10}{3},\,0,\,0,\,0,\,0,\,0,\,0,\,0,\,0,\,0,\,0,\,0,\,0,\,0,\,0,\,0\right) ,x \rangle \leq \frac{43861122710991465}{36028797018963968} \\
-\frac{40066579960122721}{36028797018963968} &\leq \langle \left(0,\,0,\,0,\,1,\,0,\,0,\,0,\,0,\,0,\,0,\,0,\,0,\,0,\,0,\,0,\,0,\,0,\,0,\,0,\,0\right) ,x \rangle \leq \frac{20040505775302443}{18014398509481984} \\
-\frac{2361185124158778945083}{2361183241434822606848} &\leq \langle \left(0,\,0,\,0,\,0,\,0,\,0,\,1,\,0,\,0,\,0,\,0,\,0,\,0,\,0,\,0,\,0,\,0,\,0,\,0,\,0\right) ,x \rangle \leq \frac{2361185124159496218833}{2361183241434822606848} \\
-\frac{2417851643085731713983019}{2417851639229258349412352} &\leq \langle \left(0,\,0,\,0,\,0,\,0,\,0,\,0,\,0,\,1,\,0,\,0,\,0,\,0,\,0,\,0,\,0,\,0,\,0,\,0,\,0\right) ,x \rangle \leq \frac{2417851643085731713983019}{2417851639229258349412352} \\
-\frac{2361185035008561037417}{2361183241434822606848} &\leq \langle \left(0,\,0,\,0,\,0,\,0,\,0,\,0,\,0,\,0,\,0,\,0,\,0,\,0,\,0,\,1,\,0,\,0,\,0,\,0,\,0\right) ,x \rangle \leq \frac{2361185035008561037417}{2361183241434822606848} \\
-\frac{76859734288972473}{72057594037927936} &\leq \langle \left(0,\,1,\,0,\,0,\,\frac{1}{1000},\,0,\,0,\,0,\,0,\,0,\,0,\,0,\,0,\,0,\,0,\,0,\,0,\,0,\,0,\,0\right) ,x \rangle \leq \frac{38648172661138361}{36028797018963968} \\
-\frac{18889473462374272866723}{18889465931478580854784} &\leq \langle \left(0,\,0,\,0,\,0,\,0,\,1,\,0,\,0,\,0,\,0,\,0,\,0,\,0,\,0,\,0,\,0,\,0,\,0,\,0,\,0\right) ,x \rangle \leq \frac{18889473462375840754951}{18889465931478580854784} \\
\end{align*}}
\subsection{Counterexample based on inequality $9$}
If $H$ is the polyhedron obtained from $N$ by removing inequality $I_{9}$, then $H$ is a bounded Hirsch polytope and $N = H \cap C$, where $C$ is the combinatorial $20$-cube defined by
{\footnotesize\begin{align*}
-\frac{4344845935841299}{1125899906842624} &\leq \langle \left(-100,\,\frac{3}{100},\,\frac{1}{50},\,0,\,-30,\,0,\,0,\,0,\,0,\,0,\,0,\,0,\,0,\,0,\,0,\,0,\,0,\,0,\,0,\,0\right) ,x \rangle \leq 1 \\
-\frac{1180593414291149733993}{1180591620717411303424} &\leq \langle \left(0,\,0,\,0,\,0,\,0,\,0,\,0,\,0,\,0,\,0,\,0,\,0,\,0,\,0,\,0,\,1,\,0,\,0,\,0,\,0\right) ,x \rangle \leq \frac{1180593414291149733993}{1180591620717411303424} \\
-\frac{2361185124158778945083}{2361183241434822606848} &\leq \langle \left(0,\,0,\,0,\,0,\,0,\,0,\,1,\,0,\,0,\,0,\,0,\,0,\,0,\,0,\,0,\,0,\,0,\,0,\,0,\,0\right) ,x \rangle \leq \frac{2361185124159496218833}{2361183241434822606848} \\
-\frac{4722372222305608191517}{4722366482869645213696} &\leq \langle \left(0,\,0,\,0,\,0,\,0,\,0,\,0,\,0,\,0,\,0,\,0,\,0,\,0,\,0,\,0,\,0,\,0,\,0,\,1,\,0\right) ,x \rangle \leq \frac{4722372222305608191517}{4722366482869645213696} \\
-\frac{18889470523027351237041}{18889465931478580854784} &\leq \langle \left(0,\,0,\,0,\,0,\,0,\,0,\,0,\,0,\,0,\,0,\,0,\,0,\,0,\,0,\,0,\,0,\,0,\,0,\,0,\,1\right) ,x \rangle \leq \frac{18889470523027351237041}{18889465931478580854784} \\
-\frac{590297603932444082281}{590295810358705651712} &\leq \langle \left(0,\,0,\,0,\,0,\,0,\,0,\,0,\,0,\,0,\,0,\,0,\,0,\,0,\,0,\,0,\,0,\,1,\,0,\,0,\,0\right) ,x \rangle \leq \frac{590297603932444082281}{590295810358705651712} \\
-\frac{9444738705175253405213}{9444732965739290427392} &\leq \langle \left(0,\,0,\,0,\,0,\,0,\,0,\,0,\,0,\,0,\,0,\,0,\,0,\,0,\,0,\,0,\,0,\,0,\,1,\,0,\,0\right) ,x \rangle \leq \frac{9444738705175253405213}{9444732965739290427392} \\
-\frac{39614081257140112583132697033}{39614081257132168796771975168} &\leq \langle \left(0,\,0,\,0,\,0,\,0,\,0,\,0,\,0,\,0,\,0,\,1,\,0,\,0,\,0,\,0,\,0,\,0,\,0,\,0,\,0\right) ,x \rangle \leq \frac{9903520314285028145783175131}{9903520314283042199192993792} \\
-\frac{2417851645399615732725419}{2417851639229258349412352} &\leq \langle \left(0,\,0,\,0,\,0,\,0,\,0,\,0,\,0,\,0,\,0,\,0,\,0,\,1,\,0,\,0,\,0,\,0,\,0,\,0,\,0\right) ,x \rangle \leq \frac{2417851645399615732725419}{2417851639229258349412352} \\
-\frac{147579558414782802697}{147573952589676412928} &\leq \langle \left(0,\,0,\,0,\,0,\,0,\,0,\,0,\,0,\,0,\,0,\,0,\,0,\,0,\,1,\,0,\,0,\,0,\,0,\,0,\,0\right) ,x \rangle \leq \frac{147579558414782802697}{147573952589676412928} \\
-\frac{9671406563106733594076309}{9671406556917033397649408} &\leq \langle \left(0,\,0,\,0,\,0,\,0,\,0,\,0,\,0,\,0,\,0,\,0,\,1,\,0,\,0,\,0,\,0,\,0,\,0,\,0,\,0\right) ,x \rangle \leq \frac{9671406563106733594076309}{9671406556917033397649408} \\
-\frac{309485009827534768921207957}{309485009821345068724781056} &\leq \langle \left(0,\,0,\,0,\,0,\,0,\,0,\,0,\,0,\,0,\,1,\,0,\,0,\,0,\,0,\,0,\,0,\,0,\,0,\,0,\,0\right) ,x \rangle \leq \frac{309485009827534768921207957}{309485009821345068724781056} \\
-\frac{43146496257989641}{36028797018963968} &\leq \langle \left(1,\,0,\,0,\,0,\,-\frac{10}{3},\,0,\,0,\,0,\,0,\,0,\,0,\,0,\,0,\,0,\,0,\,0,\,0,\,0,\,0,\,0\right) ,x \rangle \leq \frac{10973434459539691}{9007199254740992} \\
-\frac{40066579960122721}{36028797018963968} &\leq \langle \left(0,\,0,\,0,\,1,\,0,\,0,\,0,\,0,\,0,\,0,\,0,\,0,\,0,\,0,\,0,\,0,\,0,\,0,\,0,\,0\right) ,x \rangle \leq \frac{20040505775302443}{18014398509481984} \\
-\frac{1180593503442084915409}{1180591620717411303424} &\leq \langle \left(0,\,0,\,0,\,0,\,0,\,0,\,0,\,1,\,0,\,0,\,0,\,0,\,0,\,0,\,0,\,0,\,0,\,0,\,0,\,0\right) ,x \rangle \leq \frac{1180593503442084915409}{1180591620717411303424} \\
-\frac{75849975776390845}{72057594037927936} &\leq \langle \left(0,\,0,\,1,\,0,\,\frac{1}{1500},\,0,\,0,\,0,\,0,\,0,\,0,\,0,\,0,\,0,\,0,\,0,\,0,\,0,\,0,\,0\right) ,x \rangle \leq \frac{20195319152857533}{18014398509481984} \\
-\frac{2417851643085731713983019}{2417851639229258349412352} &\leq \langle \left(0,\,0,\,0,\,0,\,0,\,0,\,0,\,0,\,1,\,0,\,0,\,0,\,0,\,0,\,0,\,0,\,0,\,0,\,0,\,0\right) ,x \rangle \leq \frac{2417851643085731713983019}{2417851639229258349412352} \\
-\frac{2361185035008561037417}{2361183241434822606848} &\leq \langle \left(0,\,0,\,0,\,0,\,0,\,0,\,0,\,0,\,0,\,0,\,0,\,0,\,0,\,0,\,1,\,0,\,0,\,0,\,0,\,0\right) ,x \rangle \leq \frac{2361185035008561037417}{2361183241434822606848} \\
-\frac{76859734288972473}{72057594037927936} &\leq \langle \left(0,\,1,\,0,\,0,\,\frac{1}{1000},\,0,\,0,\,0,\,0,\,0,\,0,\,0,\,0,\,0,\,0,\,0,\,0,\,0,\,0,\,0\right) ,x \rangle \leq \frac{38648172661138361}{36028797018963968} \\
-\frac{18889473462374272866723}{18889465931478580854784} &\leq \langle \left(0,\,0,\,0,\,0,\,0,\,1,\,0,\,0,\,0,\,0,\,0,\,0,\,0,\,0,\,0,\,0,\,0,\,0,\,0,\,0\right) ,x \rangle \leq \frac{18889473462375840754951}{18889465931478580854784} \\
\end{align*}}
\subsection{Counterexample based on inequality $10$}
If $H$ is the polyhedron obtained from $N$ by removing inequality $I_{10}$, then $H$ is a bounded Hirsch polytope and $N = H \cap C$, where $C$ is the combinatorial $20$-cube defined by
{\footnotesize\begin{align*}
-\frac{76859732776487511}{72057594037927936} &\leq \langle \left(0,\,1,\,0,\,0,\,0,\,0,\,0,\,0,\,0,\,0,\,0,\,0,\,0,\,0,\,\frac{3}{2000},\,0,\,0,\,0,\,0,\,0\right) ,x \rangle \leq \frac{38648171886974793}{36028797018963968} \\
-\frac{4722372222305608191517}{4722366482869645213696} &\leq \langle \left(0,\,0,\,0,\,0,\,0,\,0,\,0,\,0,\,0,\,0,\,0,\,0,\,0,\,0,\,0,\,0,\,0,\,0,\,1,\,0\right) ,x \rangle \leq \frac{4722372222305608191517}{4722366482869645213696} \\
-\frac{2361185124158778945083}{2361183241434822606848} &\leq \langle \left(0,\,0,\,0,\,0,\,0,\,0,\,1,\,0,\,0,\,0,\,0,\,0,\,0,\,0,\,0,\,0,\,0,\,0,\,0,\,0\right) ,x \rangle \leq \frac{2361185124159496218833}{2361183241434822606848} \\
-\frac{153120913999410017}{144115188075855872} &\leq \langle \left(0,\,0,\,0,\,0,\,1,\,0,\,0,\,0,\,0,\,0,\,0,\,0,\,0,\,0,\,-\frac{159}{5},\,0,\,0,\,0,\,0,\,0\right) ,x \rangle \leq \frac{76542601517243887}{72057594037927936} \\
-\frac{9444738705175253405213}{9444732965739290427392} &\leq \langle \left(0,\,0,\,0,\,0,\,0,\,0,\,0,\,0,\,0,\,0,\,0,\,0,\,0,\,0,\,0,\,0,\,0,\,1,\,0,\,0\right) ,x \rangle \leq \frac{9444738705175253405213}{9444732965739290427392} \\
-\frac{8941986646900661}{2251799813685248} &\leq \langle \left(-100,\,\frac{3}{2000},\,-\frac{7}{2000},\,0,\,-\frac{159}{5},\,0,\,0,\,0,\,0,\,0,\,0,\,0,\,0,\,0,\,-1,\,0,\,0,\,0,\,0,\,0\right) ,x \rangle \leq 1 \\
-\frac{18889470523027351237041}{18889465931478580854784} &\leq \langle \left(0,\,0,\,0,\,0,\,0,\,0,\,0,\,0,\,0,\,0,\,0,\,0,\,0,\,0,\,0,\,0,\,0,\,0,\,0,\,1\right) ,x \rangle \leq \frac{18889470523027351237041}{18889465931478580854784} \\
-\frac{590297603932444082281}{590295810358705651712} &\leq \langle \left(0,\,0,\,0,\,0,\,0,\,0,\,0,\,0,\,0,\,0,\,0,\,0,\,0,\,0,\,0,\,0,\,1,\,0,\,0,\,0\right) ,x \rangle \leq \frac{590297603932444082281}{590295810358705651712} \\
-\frac{39614081257140112583132697033}{39614081257132168796771975168} &\leq \langle \left(0,\,0,\,0,\,0,\,0,\,0,\,0,\,0,\,0,\,0,\,1,\,0,\,0,\,0,\,0,\,0,\,0,\,0,\,0,\,0\right) ,x \rangle \leq \frac{9903520314285028145783175131}{9903520314283042199192993792} \\
-\frac{151699952049529573}{144115188075855872} &\leq \langle \left(0,\,0,\,1,\,0,\,0,\,0,\,0,\,0,\,0,\,0,\,0,\,0,\,0,\,0,\,-\frac{7}{2000},\,0,\,0,\,0,\,0,\,0\right) ,x \rangle \leq \frac{40390636902147947}{36028797018963968} \\
-\frac{2417851645399615732725419}{2417851639229258349412352} &\leq \langle \left(0,\,0,\,0,\,0,\,0,\,0,\,0,\,0,\,0,\,0,\,0,\,0,\,1,\,0,\,0,\,0,\,0,\,0,\,0,\,0\right) ,x \rangle \leq \frac{2417851645399615732725419}{2417851639229258349412352} \\
-\frac{147579558414782802697}{147573952589676412928} &\leq \langle \left(0,\,0,\,0,\,0,\,0,\,0,\,0,\,0,\,0,\,0,\,0,\,0,\,0,\,1,\,0,\,0,\,0,\,0,\,0,\,0\right) ,x \rangle \leq \frac{147579558414782802697}{147573952589676412928} \\
-\frac{9671406563106733594076309}{9671406556917033397649408} &\leq \langle \left(0,\,0,\,0,\,0,\,0,\,0,\,0,\,0,\,0,\,0,\,0,\,1,\,0,\,0,\,0,\,0,\,0,\,0,\,0,\,0\right) ,x \rangle \leq \frac{9671406563106733594076309}{9671406556917033397649408} \\
-\frac{309485009827534768921207957}{309485009821345068724781056} &\leq \langle \left(0,\,0,\,0,\,0,\,0,\,0,\,0,\,0,\,0,\,1,\,0,\,0,\,0,\,0,\,0,\,0,\,0,\,0,\,0,\,0\right) ,x \rangle \leq \frac{309485009827534768921207957}{309485009821345068724781056} \\
-\frac{1180593503442084915409}{1180591620717411303424} &\leq \langle \left(0,\,0,\,0,\,0,\,0,\,0,\,0,\,1,\,0,\,0,\,0,\,0,\,0,\,0,\,0,\,0,\,0,\,0,\,0,\,0\right) ,x \rangle \leq \frac{1180593503442084915409}{1180591620717411303424} \\
-\frac{40066579960122721}{36028797018963968} &\leq \langle \left(0,\,0,\,0,\,1,\,0,\,0,\,0,\,0,\,0,\,0,\,0,\,0,\,0,\,0,\,0,\,0,\,0,\,0,\,0,\,0\right) ,x \rangle \leq \frac{20040505775302443}{18014398509481984} \\
-\frac{2417851643085731713983019}{2417851639229258349412352} &\leq \langle \left(0,\,0,\,0,\,0,\,0,\,0,\,0,\,0,\,1,\,0,\,0,\,0,\,0,\,0,\,0,\,0,\,0,\,0,\,0,\,0\right) ,x \rangle \leq \frac{2417851643085731713983019}{2417851639229258349412352} \\
-\frac{1180593414291149733993}{1180591620717411303424} &\leq \langle \left(0,\,0,\,0,\,0,\,0,\,0,\,0,\,0,\,0,\,0,\,0,\,0,\,0,\,0,\,0,\,1,\,0,\,0,\,0,\,0\right) ,x \rangle \leq \frac{1180593414291149733993}{1180591620717411303424} \\
-\frac{291132536530816775}{288230376151711744} &\leq \langle \left(1,\,0,\,0,\,0,\,0,\,0,\,0,\,0,\,0,\,0,\,0,\,0,\,0,\,0,\,-100,\,0,\,0,\,0,\,0,\,0\right) ,x \rangle \leq \frac{582225359826457723}{576460752303423488} \\
-\frac{18889473462374272866723}{18889465931478580854784} &\leq \langle \left(0,\,0,\,0,\,0,\,0,\,1,\,0,\,0,\,0,\,0,\,0,\,0,\,0,\,0,\,0,\,0,\,0,\,0,\,0,\,0\right) ,x \rangle \leq \frac{18889473462375840754951}{18889465931478580854784} \\
\end{align*}}
\subsection{Counterexample based on inequality $11$}
If $H$ is the polyhedron obtained from $N$ by removing inequality $I_{11}$, then $H$ is a bounded Hirsch polytope and $N = H \cap C$, where $C$ is the combinatorial $20$-cube defined by
{\footnotesize\begin{align*}
-\frac{1122197785160553}{1099511627776} &\leq \langle \bigg(-100,\,-\frac{3}{2000},\,-\frac{7}{2000},\,0,\,\frac{159}{5},\,-10000000,\,-10000000,\,-10000000,\,-10000000000,\\&\,-100000000000,\,-100000000000,\,-100000000000,\,-100000000000,\,-1,\,0,\,0,\,0,\,0,\,0,\,0\bigg) ,x \rangle \leq 1 \\
-\frac{1180593414291149733993}{1180591620717411303424} &\leq \langle \left(0,\,0,\,0,\,0,\,0,\,0,\,0,\,0,\,0,\,0,\,0,\,0,\,0,\,0,\,0,\,1,\,0,\,0,\,0,\,0\right) ,x \rangle \leq \frac{1180593414291149733993}{1180591620717411303424} \\
-\frac{2361185035008561037417}{2361183241434822606848} &\leq \langle \left(0,\,0,\,0,\,0,\,0,\,0,\,0,\,0,\,0,\,0,\,0,\,0,\,0,\,0,\,1,\,0,\,0,\,0,\,0,\,0\right) ,x \rangle \leq \frac{2361185035008561037417}{2361183241434822606848} \\
-\frac{18889470523027351237041}{18889465931478580854784} &\leq \langle \left(0,\,0,\,0,\,0,\,0,\,0,\,0,\,0,\,0,\,0,\,0,\,0,\,0,\,0,\,0,\,0,\,0,\,0,\,0,\,1\right) ,x \rangle \leq \frac{18889470523027351237041}{18889465931478580854784} \\
-\frac{3350128234806005}{8796093022208} &\leq \langle \left(0,\,0,\,0,\,0,\,0,\,0,\,1,\,0,\,0,\,0,\,0,\,0,\,0,\,-10000000,\,0,\,0,\,0,\,0,\,0,\,0\right) ,x \rangle \leq \frac{1675064117393089}{4398046511104} \\
-\frac{4722372222305608191517}{4722366482869645213696} &\leq \langle \left(0,\,0,\,0,\,0,\,0,\,0,\,0,\,0,\,0,\,0,\,0,\,0,\,0,\,0,\,0,\,0,\,0,\,0,\,1,\,0\right) ,x \rangle \leq \frac{4722372222305608191517}{4722366482869645213696} \\
-\frac{3350128234803823}{8796093022208} &\leq \langle \left(0,\,0,\,0,\,0,\,0,\,0,\,0,\,1,\,0,\,0,\,0,\,0,\,0,\,-10000000,\,0,\,0,\,0,\,0,\,0,\,0\right) ,x \rangle \leq \frac{1675064117394113}{4398046511104} \\
-\frac{1631514129823731}{4294967296} &\leq \langle \left(0,\,0,\,0,\,0,\,0,\,0,\,0,\,0,\,1,\,0,\,0,\,0,\,0,\,-10000000000,\,0,\,0,\,0,\,0,\,0,\,0\right) ,x \rangle \leq \frac{3263028259623969}{8589934592} \\
-\frac{9444738705175253405213}{9444732965739290427392} &\leq \langle \left(0,\,0,\,0,\,0,\,0,\,0,\,0,\,0,\,0,\,0,\,0,\,0,\,0,\,0,\,0,\,0,\,0,\,1,\,0,\,0\right) ,x \rangle \leq \frac{9444738705175253405213}{9444732965739290427392} \\
-\frac{4078775660882911}{1073741824} &\leq \langle \left(0,\,0,\,0,\,0,\,0,\,0,\,0,\,0,\,0,\,0,\,0,\,0,\,1,\,-100000000000,\,0,\,0,\,0,\,0,\,0,\,0\right) ,x \rangle \leq \frac{8157551321707091}{2147483648} \\
-\frac{6700256469614191}{17592186044416} &\leq \langle \left(0,\,0,\,0,\,0,\,0,\,1,\,0,\,0,\,0,\,0,\,0,\,0,\,0,\,-10000000,\,0,\,0,\,0,\,0,\,0,\,0\right) ,x \rangle \leq \frac{1675064117392577}{4398046511104} \\
-\frac{4078775660882911}{1073741824} &\leq \langle \left(0,\,0,\,0,\,0,\,0,\,0,\,0,\,0,\,0,\,0,\,0,\,1,\,0,\,-100000000000,\,0,\,0,\,0,\,0,\,0,\,0\right) ,x \rangle \leq \frac{8157551321707091}{2147483648} \\
-\frac{4078775660882911}{1073741824} &\leq \langle \left(0,\,0,\,0,\,0,\,0,\,0,\,0,\,0,\,0,\,1,\,0,\,0,\,0,\,-100000000000,\,0,\,0,\,0,\,0,\,0,\,0\right) ,x \rangle \leq \frac{8157551321707091}{2147483648} \\
-\frac{292205526492713759}{288230376151711744} &\leq \langle \left(1,\,0,\,0,\,0,\,0,\,0,\,0,\,0,\,0,\,0,\,0,\,0,\,0,\,-100,\,0,\,0,\,0,\,0,\,0,\,0\right) ,x \rangle \leq \frac{582225360191781879}{576460752303423488} \\
-\frac{590297603932444082281}{590295810358705651712} &\leq \langle \left(0,\,0,\,0,\,0,\,0,\,0,\,0,\,0,\,0,\,0,\,0,\,0,\,0,\,0,\,0,\,0,\,1,\,0,\,0,\,0\right) ,x \rangle \leq \frac{590297603932444082281}{590295810358705651712} \\
-\frac{18962495828811293}{18014398509481984} &\leq \langle \left(0,\,0,\,1,\,0,\,0,\,0,\,0,\,0,\,0,\,0,\,0,\,0,\,0,\,-\frac{7}{2000},\,0,\,0,\,0,\,0,\,0,\,0\right) ,x \rangle \leq \frac{40390640916459081}{36028797018963968} \\
-\frac{4078775660882911}{1073741824} &\leq \langle \left(0,\,0,\,0,\,0,\,0,\,0,\,0,\,0,\,0,\,0,\,1,\,0,\,0,\,-100000000000,\,0,\,0,\,0,\,0,\,0,\,0\right) ,x \rangle \leq \frac{8157551321707091}{2147483648} \\
-\frac{19140120201425635}{18014398509481984} &\leq \langle \left(0,\,0,\,0,\,0,\,1,\,0,\,0,\,0,\,0,\,0,\,0,\,0,\,0,\,\frac{159}{5},\,0,\,0,\,0,\,0,\,0,\,0\right) ,x \rangle \leq \frac{76542601517243887}{72057594037927936} \\
-\frac{9607467100024973}{9007199254740992} &\leq \langle \left(0,\,1,\,0,\,0,\,0,\,0,\,0,\,0,\,0,\,0,\,0,\,0,\,0,\,-\frac{3}{2000},\,0,\,0,\,0,\,0,\,0,\,0\right) ,x \rangle \leq \frac{9662042973792413}{9007199254740992} \\
-\frac{40066579960122721}{36028797018963968} &\leq \langle \left(0,\,0,\,0,\,1,\,0,\,0,\,0,\,0,\,0,\,0,\,0,\,0,\,0,\,0,\,0,\,0,\,0,\,0,\,0,\,0\right) ,x \rangle \leq \frac{20040505775302443}{18014398509481984} \\
\end{align*}}
\subsection{Counterexample based on inequality $20$}
If $H$ is the polyhedron obtained from $N$ by removing inequality $I_{20}$, then $H$ is a bounded Hirsch polytope and $N = H \cap C$, where $C$ is the combinatorial $20$-cube defined by
{\footnotesize\begin{align*}
-\frac{151699951751958945}{144115188075855872} &\leq \langle \left(0,\,0,\,1,\,0,\,0,\,0,\,0,\,0,\,0,\,0,\,0,\,0,\,0,\,0,\,0,\,0,\,0,\,0,\,0,\,0\right) ,x \rangle \leq \frac{80781273703537389}{72057594037927936} \\
-\frac{4722372222305608191517}{4722366482869645213696} &\leq \langle \left(0,\,0,\,0,\,0,\,0,\,0,\,0,\,0,\,0,\,0,\,0,\,0,\,0,\,0,\,0,\,0,\,0,\,0,\,1,\,0\right) ,x \rangle \leq \frac{4722372222305608191517}{4722366482869645213696} \\
-\frac{1180593414291149733993}{1180591620717411303424} &\leq \langle \left(0,\,0,\,0,\,0,\,0,\,0,\,0,\,0,\,0,\,0,\,0,\,0,\,0,\,0,\,0,\,1,\,0,\,0,\,0,\,0\right) ,x \rangle \leq \frac{1180593414291149733993}{1180591620717411303424} \\
-\frac{2361185124158778945083}{2361183241434822606848} &\leq \langle \left(0,\,0,\,0,\,0,\,0,\,0,\,1,\,0,\,0,\,0,\,0,\,0,\,0,\,0,\,0,\,0,\,0,\,0,\,0,\,0\right) ,x \rangle \leq \frac{2361185124159496218833}{2361183241434822606848} \\
-\frac{38429866348206703}{36028797018963968} &\leq \langle \left(0,\,1,\,0,\,0,\,0,\,0,\,0,\,0,\,-30,\,0,\,0,\,0,\,0,\,0,\,0,\,0,\,0,\,0,\,0,\,0\right) ,x \rangle \leq \frac{77296343773103709}{72057594037927936} \\
-\frac{18889470523027351237041}{18889465931478580854784} &\leq \langle \left(0,\,0,\,0,\,0,\,0,\,0,\,0,\,0,\,0,\,0,\,0,\,0,\,0,\,0,\,0,\,0,\,0,\,0,\,0,\,1\right) ,x \rangle \leq \frac{18889470523027351237041}{18889465931478580854784} \\
-\frac{590297603932444082281}{590295810358705651712} &\leq \langle \left(0,\,0,\,0,\,0,\,0,\,0,\,0,\,0,\,0,\,0,\,0,\,0,\,0,\,0,\,0,\,0,\,1,\,0,\,0,\,0\right) ,x \rangle \leq \frac{590297603932444082281}{590295810358705651712} \\
-\frac{9444738705175253405213}{9444732965739290427392} &\leq \langle \left(0,\,0,\,0,\,0,\,0,\,0,\,0,\,0,\,0,\,0,\,0,\,0,\,0,\,0,\,0,\,0,\,0,\,1,\,0,\,0\right) ,x \rangle \leq \frac{9444738705175253405213}{9444732965739290427392} \\
-\frac{39614081257140112583132697033}{39614081257132168796771975168} &\leq \langle \left(0,\,0,\,0,\,0,\,0,\,0,\,0,\,0,\,0,\,0,\,1,\,0,\,0,\,0,\,0,\,0,\,0,\,0,\,0,\,0\right) ,x \rangle \leq \frac{9903520314285028145783175131}{9903520314283042199192993792} \\
-\frac{2417851645399615732725419}{2417851639229258349412352} &\leq \langle \left(0,\,0,\,0,\,0,\,0,\,0,\,0,\,0,\,0,\,0,\,0,\,0,\,1,\,0,\,0,\,0,\,0,\,0,\,0,\,0\right) ,x \rangle \leq \frac{2417851645399615732725419}{2417851639229258349412352} \\
-\frac{147579558414782802697}{147573952589676412928} &\leq \langle \left(0,\,0,\,0,\,0,\,0,\,0,\,0,\,0,\,0,\,0,\,0,\,0,\,0,\,1,\,0,\,0,\,0,\,0,\,0,\,0\right) ,x \rangle \leq \frac{147579558414782802697}{147573952589676412928} \\
-\frac{1176672303808799}{281474976710656} &\leq \langle \left(100,\,-30,\,0,\,0,\,0,\,0,\,0,\,-1,\,-1,\,0,\,0,\,0,\,0,\,0,\,0,\,0,\,0,\,0,\,0,\,0\right) ,x \rangle \leq 1 \\
-\frac{9671406563106733594076309}{9671406556917033397649408} &\leq \langle \left(0,\,0,\,0,\,0,\,0,\,0,\,0,\,0,\,0,\,0,\,0,\,1,\,0,\,0,\,0,\,0,\,0,\,0,\,0,\,0\right) ,x \rangle \leq \frac{9671406563106733594076309}{9671406556917033397649408} \\
-\frac{309485009827534768921207957}{309485009821345068724781056} &\leq \langle \left(0,\,0,\,0,\,0,\,0,\,0,\,0,\,0,\,0,\,1,\,0,\,0,\,0,\,0,\,0,\,0,\,0,\,0,\,0,\,0\right) ,x \rangle \leq \frac{309485009827534768921207957}{309485009821345068724781056} \\
-\frac{1180593503442084915409}{1180591620717411303424} &\leq \langle \left(0,\,0,\,0,\,0,\,0,\,0,\,0,\,1,\,-1,\,0,\,0,\,0,\,0,\,0,\,0,\,0,\,0,\,0,\,0,\,0\right) ,x \rangle \leq \frac{1180593507207534262633}{1180591620717411303424} \\
-\frac{40066579960122721}{36028797018963968} &\leq \langle \left(0,\,0,\,0,\,1,\,0,\,0,\,0,\,0,\,0,\,0,\,0,\,0,\,0,\,0,\,0,\,0,\,0,\,0,\,0,\,0\right) ,x \rangle \leq \frac{20040505775302443}{18014398509481984} \\
-\frac{18889473462374272866723}{18889465931478580854784} &\leq \langle \left(0,\,0,\,0,\,0,\,0,\,1,\,0,\,0,\,0,\,0,\,0,\,0,\,0,\,0,\,0,\,0,\,0,\,0,\,0,\,0\right) ,x \rangle \leq \frac{18889473462375840754951}{18889465931478580854784} \\
-\frac{2361185035008561037417}{2361183241434822606848} &\leq \langle \left(0,\,0,\,0,\,0,\,0,\,0,\,0,\,0,\,0,\,0,\,0,\,0,\,0,\,0,\,1,\,0,\,0,\,0,\,0,\,0\right) ,x \rangle \leq \frac{2361185035008561037417}{2361183241434822606848} \\
-\frac{582225365103857101}{576460752303423488} &\leq \langle \left(1,\,0,\,0,\,0,\,0,\,0,\,0,\,0,\,100,\,0,\,0,\,0,\,0,\,0,\,0,\,0,\,0,\,0,\,0,\,0\right) ,x \rangle \leq \frac{291112726061874127}{288230376151711744} \\
-\frac{2392514266058385}{2251799813685248} &\leq \langle \left(0,\,0,\,0,\,0,\,1,\,0,\,0,\,0,\,0,\,0,\,0,\,0,\,0,\,0,\,0,\,0,\,0,\,0,\,0,\,0\right) ,x \rangle \leq \frac{76542601517243887}{72057594037927936} \\
\end{align*}}
\subsection{Counterexample based on inequality $21$}
If $H$ is the polyhedron obtained from $N$ by removing inequality $I_{21}$, then $H$ is a bounded Hirsch polytope and $N = H \cap C$, where $C$ is the combinatorial $20$-cube defined by
{\footnotesize\begin{align*}
-\frac{6837637516946527}{2251799813685248} &\leq \langle \left(100,\,-4,\,15,\,0,\,0,\,0,\,0,\,0,\,-1,\,-1,\,0,\,0,\,0,\,0,\,0,\,0,\,0,\,0,\,0,\,0\right) ,x \rangle \leq 1 \\
-\frac{76859732694389839}{72057594037927936} &\leq \langle \left(0,\,1,\,0,\,0,\,0,\,0,\,0,\,0,\,0,\,-4,\,0,\,0,\,0,\,0,\,0,\,0,\,0,\,0,\,0,\,0\right) ,x \rangle \leq \frac{77296343773162177}{72057594037927936} \\
-\frac{2392514266058385}{2251799813685248} &\leq \langle \left(0,\,0,\,0,\,0,\,1,\,0,\,0,\,0,\,0,\,0,\,0,\,0,\,0,\,0,\,0,\,0,\,0,\,0,\,0,\,0\right) ,x \rangle \leq \frac{76542601517243887}{72057594037927936} \\
-\frac{2361185124158778945083}{2361183241434822606848} &\leq \langle \left(0,\,0,\,0,\,0,\,0,\,0,\,1,\,0,\,0,\,0,\,0,\,0,\,0,\,0,\,0,\,0,\,0,\,0,\,0,\,0\right) ,x \rangle \leq \frac{2361185124159496218833}{2361183241434822606848} \\
-\frac{4722372222305608191517}{4722366482869645213696} &\leq \langle \left(0,\,0,\,0,\,0,\,0,\,0,\,0,\,0,\,0,\,0,\,0,\,0,\,0,\,0,\,0,\,0,\,0,\,0,\,1,\,0\right) ,x \rangle \leq \frac{4722372222305608191517}{4722366482869645213696} \\
-\frac{18889470523027351237041}{18889465931478580854784} &\leq \langle \left(0,\,0,\,0,\,0,\,0,\,0,\,0,\,0,\,0,\,0,\,0,\,0,\,0,\,0,\,0,\,0,\,0,\,0,\,0,\,1\right) ,x \rangle \leq \frac{18889470523027351237041}{18889465931478580854784} \\
-\frac{72778170637900421}{72057594037927936} &\leq \langle \left(1,\,0,\,0,\,0,\,0,\,0,\,0,\,0,\,0,\,100,\,0,\,0,\,0,\,0,\,0,\,0,\,0,\,0,\,0,\,0\right) ,x \rangle \leq \frac{582225378290904667}{576460752303423488} \\
-\frac{9444738705175253405213}{9444732965739290427392} &\leq \langle \left(0,\,0,\,0,\,0,\,0,\,0,\,0,\,0,\,0,\,0,\,0,\,0,\,0,\,0,\,0,\,0,\,0,\,1,\,0,\,0\right) ,x \rangle \leq \frac{9444738705175253405213}{9444732965739290427392} \\
-\frac{75849975875943419}{72057594037927936} &\leq \langle \left(0,\,0,\,1,\,0,\,0,\,0,\,0,\,0,\,0,\,15,\,0,\,0,\,0,\,0,\,0,\,0,\,0,\,0,\,0,\,0\right) ,x \rangle \leq \frac{40390636852017309}{36028797018963968} \\
-\frac{39614081257140112583132697033}{39614081257132168796771975168} &\leq \langle \left(0,\,0,\,0,\,0,\,0,\,0,\,0,\,0,\,0,\,0,\,1,\,0,\,0,\,0,\,0,\,0,\,0,\,0,\,0,\,0\right) ,x \rangle \leq \frac{9903520314285028145783175131}{9903520314283042199192993792} \\
-\frac{2417851645399615732725419}{2417851639229258349412352} &\leq \langle \left(0,\,0,\,0,\,0,\,0,\,0,\,0,\,0,\,0,\,0,\,0,\,0,\,1,\,0,\,0,\,0,\,0,\,0,\,0,\,0\right) ,x \rangle \leq \frac{2417851645399615732725419}{2417851639229258349412352} \\
-\frac{147579558414782802697}{147573952589676412928} &\leq \langle \left(0,\,0,\,0,\,0,\,0,\,0,\,0,\,0,\,0,\,0,\,0,\,0,\,0,\,1,\,0,\,0,\,0,\,0,\,0,\,0\right) ,x \rangle \leq \frac{147579558414782802697}{147573952589676412928} \\
-\frac{2417851643085731713983019}{2417851639229258349412352} &\leq \langle \left(0,\,0,\,0,\,0,\,0,\,0,\,0,\,0,\,1,\,-1,\,0,\,0,\,0,\,0,\,0,\,0,\,0,\,0,\,0,\,0\right) ,x \rangle \leq \frac{9444732983816509323817}{9444732965739290427392} \\
-\frac{9671406563106733594076309}{9671406556917033397649408} &\leq \langle \left(0,\,0,\,0,\,0,\,0,\,0,\,0,\,0,\,0,\,0,\,0,\,1,\,0,\,0,\,0,\,0,\,0,\,0,\,0,\,0\right) ,x \rangle \leq \frac{9671406563106733594076309}{9671406556917033397649408} \\
-\frac{590297603932444082281}{590295810358705651712} &\leq \langle \left(0,\,0,\,0,\,0,\,0,\,0,\,0,\,0,\,0,\,0,\,0,\,0,\,0,\,0,\,0,\,0,\,1,\,0,\,0,\,0\right) ,x \rangle \leq \frac{590297603932444082281}{590295810358705651712} \\
-\frac{1180593414291149733993}{1180591620717411303424} &\leq \langle \left(0,\,0,\,0,\,0,\,0,\,0,\,0,\,0,\,0,\,0,\,0,\,0,\,0,\,0,\,0,\,1,\,0,\,0,\,0,\,0\right) ,x \rangle \leq \frac{1180593414291149733993}{1180591620717411303424} \\
-\frac{1180593503442084915409}{1180591620717411303424} &\leq \langle \left(0,\,0,\,0,\,0,\,0,\,0,\,0,\,1,\,0,\,0,\,0,\,0,\,0,\,0,\,0,\,0,\,0,\,0,\,0,\,0\right) ,x \rangle \leq \frac{1180593503442084915409}{1180591620717411303424} \\
-\frac{40066579960122721}{36028797018963968} &\leq \langle \left(0,\,0,\,0,\,1,\,0,\,0,\,0,\,0,\,0,\,0,\,0,\,0,\,0,\,0,\,0,\,0,\,0,\,0,\,0,\,0\right) ,x \rangle \leq \frac{20040505775302443}{18014398509481984} \\
-\frac{2361185035008561037417}{2361183241434822606848} &\leq \langle \left(0,\,0,\,0,\,0,\,0,\,0,\,0,\,0,\,0,\,0,\,0,\,0,\,0,\,0,\,1,\,0,\,0,\,0,\,0,\,0\right) ,x \rangle \leq \frac{2361185035008561037417}{2361183241434822606848} \\
-\frac{18889473462374272866723}{18889465931478580854784} &\leq \langle \left(0,\,0,\,0,\,0,\,0,\,1,\,0,\,0,\,0,\,0,\,0,\,0,\,0,\,0,\,0,\,0,\,0,\,0,\,0,\,0\right) ,x \rangle \leq \frac{18889473462375840754951}{18889465931478580854784} \\
\end{align*}}
\subsection{Counterexample based on inequality $22$}
If $H$ is the polyhedron obtained from $N$ by removing inequality $I_{22}$, then $H$ is a bounded Hirsch polytope and $N = H \cap C$, where $C$ is the combinatorial $20$-cube defined by
{\footnotesize\begin{align*}
-\frac{1180593414291149733993}{1180591620717411303424} &\leq \langle \left(0,\,0,\,0,\,0,\,0,\,0,\,0,\,0,\,0,\,0,\,0,\,0,\,0,\,0,\,0,\,1,\,0,\,0,\,0,\,0\right) ,x \rangle \leq \frac{1180593414291149733993}{1180591620717411303424} \\
-\frac{2361185124158778945083}{2361183241434822606848} &\leq \langle \left(0,\,0,\,0,\,0,\,0,\,0,\,1,\,0,\,0,\,0,\,0,\,0,\,0,\,0,\,0,\,0,\,0,\,0,\,0,\,0\right) ,x \rangle \leq \frac{2361185124159496218833}{2361183241434822606848} \\
-\frac{4722372222305608191517}{4722366482869645213696} &\leq \langle \left(0,\,0,\,0,\,0,\,0,\,0,\,0,\,0,\,0,\,0,\,0,\,0,\,0,\,0,\,0,\,0,\,0,\,0,\,1,\,0\right) ,x \rangle \leq \frac{4722372222305608191517}{4722366482869645213696} \\
-\frac{76859732694334399}{72057594037927936} &\leq \langle \left(0,\,1,\,0,\,0,\,0,\,0,\,0,\,0,\,0,\,0,\,0,\,0,\,0,\,0,\,0,\,0,\,0,\,0,\,0,\,0\right) ,x \rangle \leq \frac{77296343773163779}{72057594037927936} \\
-\frac{18889470523027351237041}{18889465931478580854784} &\leq \langle \left(0,\,0,\,0,\,0,\,0,\,0,\,0,\,0,\,0,\,0,\,0,\,0,\,0,\,0,\,0,\,0,\,0,\,0,\,0,\,1\right) ,x \rangle \leq \frac{18889470523027351237041}{18889465931478580854784} \\
-\frac{9444738705175253405213}{9444732965739290427392} &\leq \langle \left(0,\,0,\,0,\,0,\,0,\,0,\,0,\,0,\,0,\,0,\,0,\,0,\,0,\,0,\,0,\,0,\,0,\,1,\,0,\,0\right) ,x \rangle \leq \frac{9444738705175253405213}{9444732965739290427392} \\
-\frac{9671406563106733594076309}{9671406556917033397649408} &\leq \langle \left(0,\,0,\,0,\,0,\,0,\,0,\,0,\,0,\,0,\,0,\,0,\,1,\,0,\,0,\,0,\,0,\,0,\,0,\,0,\,0\right) ,x \rangle \leq \frac{9671406563106733594076309}{9671406556917033397649408} \\
-\frac{2417851645399615732725419}{2417851639229258349412352} &\leq \langle \left(0,\,0,\,0,\,0,\,0,\,0,\,0,\,0,\,0,\,0,\,0,\,0,\,1,\,0,\,0,\,0,\,0,\,0,\,0,\,0\right) ,x \rangle \leq \frac{2417851645399615732725419}{2417851639229258349412352} \\
-\frac{147579558414782802697}{147573952589676412928} &\leq \langle \left(0,\,0,\,0,\,0,\,0,\,0,\,0,\,0,\,0,\,0,\,0,\,0,\,0,\,1,\,0,\,0,\,0,\,0,\,0,\,0\right) ,x \rangle \leq \frac{147579558414782802697}{147573952589676412928} \\
-\frac{582225365103366801}{576460752303423488} &\leq \langle \left(1,\,0,\,0,\,0,\,0,\,0,\,0,\,0,\,0,\,0,\,100,\,0,\,0,\,0,\,0,\,0,\,0,\,0,\,0,\,0\right) ,x \rangle \leq \frac{582225396749359173}{576460752303423488} \\
-\frac{3229450158965649}{1125899906842624} &\leq \langle \left(100,\,0,\,\frac{33}{2},\,0,\,0,\,0,\,0,\,0,\,0,\,-1,\,-1,\,0,\,0,\,0,\,0,\,0,\,0,\,0,\,0,\,0\right) ,x \rangle \leq 1 \\
-\frac{590297603932444082281}{590295810358705651712} &\leq \langle \left(0,\,0,\,0,\,0,\,0,\,0,\,0,\,0,\,0,\,0,\,0,\,0,\,0,\,0,\,0,\,0,\,1,\,0,\,0,\,0\right) ,x \rangle \leq \frac{590297603932444082281}{590295810358705651712} \\
-\frac{309485009827534768921207957}{309485009821345068724781056} &\leq \langle \left(0,\,0,\,0,\,0,\,0,\,0,\,0,\,0,\,0,\,1,\,-1,\,0,\,0,\,0,\,0,\,0,\,0,\,0,\,0,\,0\right) ,x \rangle \leq \frac{604462910097456784060599}{604462909807314587353088} \\
-\frac{40066579960122721}{36028797018963968} &\leq \langle \left(0,\,0,\,0,\,1,\,0,\,0,\,0,\,0,\,0,\,0,\,0,\,0,\,0,\,0,\,0,\,0,\,0,\,0,\,0,\,0\right) ,x \rangle \leq \frac{20040505775302443}{18014398509481984} \\
-\frac{1180593503442084915409}{1180591620717411303424} &\leq \langle \left(0,\,0,\,0,\,0,\,0,\,0,\,0,\,1,\,0,\,0,\,0,\,0,\,0,\,0,\,0,\,0,\,0,\,0,\,0,\,0\right) ,x \rangle \leq \frac{1180593503442084915409}{1180591620717411303424} \\
-\frac{18889473462374272866723}{18889465931478580854784} &\leq \langle \left(0,\,0,\,0,\,0,\,0,\,1,\,0,\,0,\,0,\,0,\,0,\,0,\,0,\,0,\,0,\,0,\,0,\,0,\,0,\,0\right) ,x \rangle \leq \frac{18889473462375840754951}{18889465931478580854784} \\
-\frac{2417851643085731713983019}{2417851639229258349412352} &\leq \langle \left(0,\,0,\,0,\,0,\,0,\,0,\,0,\,0,\,1,\,0,\,0,\,0,\,0,\,0,\,0,\,0,\,0,\,0,\,0,\,0\right) ,x \rangle \leq \frac{2417851643085731713983019}{2417851639229258349412352} \\
-\frac{2361185035008561037417}{2361183241434822606848} &\leq \langle \left(0,\,0,\,0,\,0,\,0,\,0,\,0,\,0,\,0,\,0,\,0,\,0,\,0,\,0,\,1,\,0,\,0,\,0,\,0,\,0\right) ,x \rangle \leq \frac{2361185035008561037417}{2361183241434822606848} \\
-\frac{151699951751800311}{144115188075855872} &\leq \langle \left(0,\,0,\,1,\,0,\,0,\,0,\,0,\,0,\,0,\,0,\,\frac{33}{2},\,0,\,0,\,0,\,0,\,0,\,0,\,0,\,0,\,0\right) ,x \rangle \leq \frac{20195318426156259}{18014398509481984} \\
-\frac{2392514266058385}{2251799813685248} &\leq \langle \left(0,\,0,\,0,\,0,\,1,\,0,\,0,\,0,\,0,\,0,\,0,\,0,\,0,\,0,\,0,\,0,\,0,\,0,\,0,\,0\right) ,x \rangle \leq \frac{76542601517243887}{72057594037927936} \\
\end{align*}}
\subsection{Counterexample based on inequality $23$}
If $H$ is the polyhedron obtained from $N$ by removing inequality $I_{23}$, then $H$ is a bounded Hirsch polytope and $N = H \cap C$, where $C$ is the combinatorial $20$-cube defined by
{\footnotesize\begin{align*}
-\frac{39614081257140112583132711257}{39614081257132168796771975168} &\leq \langle \left(0,\,0,\,0,\,0,\,0,\,0,\,0,\,0,\,0,\,0,\,1,\,-1,\,0,\,0,\,0,\,0,\,0,\,0,\,0,\,0\right) ,x \rangle \leq \frac{9444732983816509323817}{9444732965739290427392} \\
-\frac{1180593414291149733993}{1180591620717411303424} &\leq \langle \left(0,\,0,\,0,\,0,\,0,\,0,\,0,\,0,\,0,\,0,\,0,\,0,\,0,\,0,\,0,\,1,\,0,\,0,\,0,\,0\right) ,x \rangle \leq \frac{1180593414291149733993}{1180591620717411303424} \\
-\frac{2361185124158778945083}{2361183241434822606848} &\leq \langle \left(0,\,0,\,0,\,0,\,0,\,0,\,1,\,0,\,0,\,0,\,0,\,0,\,0,\,0,\,0,\,0,\,0,\,0,\,0,\,0\right) ,x \rangle \leq \frac{2361185124159496218833}{2361183241434822606848} \\
-\frac{4722372222305608191517}{4722366482869645213696} &\leq \langle \left(0,\,0,\,0,\,0,\,0,\,0,\,0,\,0,\,0,\,0,\,0,\,0,\,0,\,0,\,0,\,0,\,0,\,0,\,1,\,0\right) ,x \rangle \leq \frac{4722372222305608191517}{4722366482869645213696} \\
-\frac{590297603932444082281}{590295810358705651712} &\leq \langle \left(0,\,0,\,0,\,0,\,0,\,0,\,0,\,0,\,0,\,0,\,0,\,0,\,0,\,0,\,0,\,0,\,1,\,0,\,0,\,0\right) ,x \rangle \leq \frac{590297603932444082281}{590295810358705651712} \\
-\frac{18889470523027351237041}{18889465931478580854784} &\leq \langle \left(0,\,0,\,0,\,0,\,0,\,0,\,0,\,0,\,0,\,0,\,0,\,0,\,0,\,0,\,0,\,0,\,0,\,0,\,0,\,1\right) ,x \rangle \leq \frac{18889470523027351237041}{18889465931478580854784} \\
-\frac{145556341275923417}{144115188075855872} &\leq \langle \left(1,\,0,\,0,\,0,\,0,\,0,\,0,\,0,\,0,\,0,\,0,\,100,\,0,\,0,\,0,\,0,\,0,\,0,\,0,\,0\right) ,x \rangle \leq \frac{291112716833134093}{288230376151711744} \\
-\frac{9444738705175253405213}{9444732965739290427392} &\leq \langle \left(0,\,0,\,0,\,0,\,0,\,0,\,0,\,0,\,0,\,0,\,0,\,0,\,0,\,0,\,0,\,0,\,0,\,1,\,0,\,0\right) ,x \rangle \leq \frac{9444738705175253405213}{9444732965739290427392} \\
-\frac{76859732694389839}{72057594037927936} &\leq \langle \left(0,\,1,\,0,\,0,\,0,\,0,\,0,\,0,\,0,\,0,\,0,\,1,\,0,\,0,\,0,\,0,\,0,\,0,\,0,\,0\right) ,x \rangle \leq \frac{38648171886581059}{36028797018963968} \\
-\frac{2417851645399615732725419}{2417851639229258349412352} &\leq \langle \left(0,\,0,\,0,\,0,\,0,\,0,\,0,\,0,\,0,\,0,\,0,\,0,\,1,\,0,\,0,\,0,\,0,\,0,\,0,\,0\right) ,x \rangle \leq \frac{2417851645399615732725419}{2417851639229258349412352} \\
-\frac{2361185035008561037417}{2361183241434822606848} &\leq \langle \left(0,\,0,\,0,\,0,\,0,\,0,\,0,\,0,\,0,\,0,\,0,\,0,\,0,\,0,\,1,\,0,\,0,\,0,\,0,\,0\right) ,x \rangle \leq \frac{2361185035008561037417}{2361183241434822606848} \\
-\frac{147579558414782802697}{147573952589676412928} &\leq \langle \left(0,\,0,\,0,\,0,\,0,\,0,\,0,\,0,\,0,\,0,\,0,\,0,\,0,\,1,\,0,\,0,\,0,\,0,\,0,\,0\right) ,x \rangle \leq \frac{147579558414782802697}{147573952589676412928} \\
-\frac{309485009827534768921207957}{309485009821345068724781056} &\leq \langle \left(0,\,0,\,0,\,0,\,0,\,0,\,0,\,0,\,0,\,1,\,0,\,0,\,0,\,0,\,0,\,0,\,0,\,0,\,0,\,0\right) ,x \rangle \leq \frac{309485009827534768921207957}{309485009821345068724781056} \\
-\frac{6512048252496539}{2251799813685248} &\leq \langle \left(100,\,1,\,16,\,0,\,0,\,0,\,0,\,0,\,0,\,0,\,-1,\,-1,\,0,\,0,\,0,\,0,\,0,\,0,\,0,\,0\right) ,x \rangle \leq 1 \\
-\frac{40066579960122721}{36028797018963968} &\leq \langle \left(0,\,0,\,0,\,1,\,0,\,0,\,0,\,0,\,0,\,0,\,0,\,0,\,0,\,0,\,0,\,0,\,0,\,0,\,0,\,0\right) ,x \rangle \leq \frac{20040505775302443}{18014398509481984} \\
-\frac{1180593503442084915409}{1180591620717411303424} &\leq \langle \left(0,\,0,\,0,\,0,\,0,\,0,\,0,\,1,\,0,\,0,\,0,\,0,\,0,\,0,\,0,\,0,\,0,\,0,\,0,\,0\right) ,x \rangle \leq \frac{1180593503442084915409}{1180591620717411303424} \\
-\frac{2392514266058385}{2251799813685248} &\leq \langle \left(0,\,0,\,0,\,0,\,1,\,0,\,0,\,0,\,0,\,0,\,0,\,0,\,0,\,0,\,0,\,0,\,0,\,0,\,0,\,0\right) ,x \rangle \leq \frac{76542601517243887}{72057594037927936} \\
-\frac{2417851643085731713983019}{2417851639229258349412352} &\leq \langle \left(0,\,0,\,0,\,0,\,0,\,0,\,0,\,0,\,1,\,0,\,0,\,0,\,0,\,0,\,0,\,0,\,0,\,0,\,0,\,0\right) ,x \rangle \leq \frac{2417851643085731713983019}{2417851639229258349412352} \\
-\frac{151699560131112911}{144115188075855872} &\leq \langle \left(0,\,0,\,1,\,0,\,0,\,0,\,0,\,0,\,0,\,0,\,0,\,16,\,0,\,0,\,0,\,0,\,0,\,0,\,0,\,0\right) ,x \rangle \leq \frac{20195318426410467}{18014398509481984} \\
-\frac{18889473462374272866723}{18889465931478580854784} &\leq \langle \left(0,\,0,\,0,\,0,\,0,\,1,\,0,\,0,\,0,\,0,\,0,\,0,\,0,\,0,\,0,\,0,\,0,\,0,\,0,\,0\right) ,x \rangle \leq \frac{18889473462375840754951}{18889465931478580854784} \\
\end{align*}}
\subsection{Counterexample based on inequality $24$}
If $H$ is the polyhedron obtained from $N$ by removing inequality $I_{24}$, then $H$ is a bounded Hirsch polytope and $N = H \cap C$, where $C$ is the combinatorial $20$-cube defined by
{\footnotesize\begin{align*}
-\frac{151699951751958945}{144115188075855872} &\leq \langle \left(0,\,0,\,1,\,0,\,0,\,0,\,0,\,0,\,0,\,0,\,0,\,0,\,0,\,0,\,0,\,0,\,0,\,0,\,0,\,0\right) ,x \rangle \leq \frac{80781273703537389}{72057594037927936} \\
-\frac{1180593414291149733993}{1180591620717411303424} &\leq \langle \left(0,\,0,\,0,\,0,\,0,\,0,\,0,\,0,\,0,\,0,\,0,\,0,\,0,\,0,\,0,\,1,\,0,\,0,\,0,\,0\right) ,x \rangle \leq \frac{1180593414291149733993}{1180591620717411303424} \\
-\frac{2361185124158778945083}{2361183241434822606848} &\leq \langle \left(0,\,0,\,0,\,0,\,0,\,0,\,1,\,0,\,0,\,0,\,0,\,0,\,0,\,0,\,0,\,0,\,0,\,0,\,0,\,0\right) ,x \rangle \leq \frac{2361185124159496218833}{2361183241434822606848} \\
-\frac{4722372222305608191517}{4722366482869645213696} &\leq \langle \left(0,\,0,\,0,\,0,\,0,\,0,\,0,\,0,\,0,\,0,\,0,\,0,\,0,\,0,\,0,\,0,\,0,\,0,\,1,\,0\right) ,x \rangle \leq \frac{4722372222305608191517}{4722366482869645213696} \\
-\frac{18889470523027351237041}{18889465931478580854784} &\leq \langle \left(0,\,0,\,0,\,0,\,0,\,0,\,0,\,0,\,0,\,0,\,0,\,0,\,0,\,0,\,0,\,0,\,0,\,0,\,0,\,1\right) ,x \rangle \leq \frac{18889470523027351237041}{18889465931478580854784} \\
-\frac{590297603932444082281}{590295810358705651712} &\leq \langle \left(0,\,0,\,0,\,0,\,0,\,0,\,0,\,0,\,0,\,0,\,0,\,0,\,0,\,0,\,0,\,0,\,1,\,0,\,0,\,0\right) ,x \rangle \leq \frac{590297603932444082281}{590295810358705651712} \\
-\frac{9444738705175253405213}{9444732965739290427392} &\leq \langle \left(0,\,0,\,0,\,0,\,0,\,0,\,0,\,0,\,0,\,0,\,0,\,0,\,0,\,0,\,0,\,0,\,0,\,1,\,0,\,0\right) ,x \rangle \leq \frac{9444738705175253405213}{9444732965739290427392} \\
-\frac{39614081257140112583132697033}{39614081257132168796771975168} &\leq \langle \left(0,\,0,\,0,\,0,\,0,\,0,\,0,\,0,\,0,\,0,\,1,\,0,\,0,\,0,\,0,\,0,\,0,\,0,\,0,\,0\right) ,x \rangle \leq \frac{9903520314285028145783175131}{9903520314283042199192993792} \\
-\frac{76859732697383609}{72057594037927936} &\leq \langle \left(0,\,1,\,0,\,0,\,0,\,0,\,0,\,0,\,0,\,0,\,0,\,0,\,\frac{55}{2},\,0,\,0,\,0,\,0,\,0,\,0,\,0\right) ,x \rangle \leq \frac{38648171886537751}{36028797018963968} \\
-\frac{147579558414782802697}{147573952589676412928} &\leq \langle \left(0,\,0,\,0,\,0,\,0,\,0,\,0,\,0,\,0,\,0,\,0,\,0,\,0,\,1,\,0,\,0,\,0,\,0,\,0,\,0\right) ,x \rangle \leq \frac{147579558414782802697}{147573952589676412928} \\
-\frac{2157210724287645}{562949953421312} &\leq \langle \left(100,\,\frac{55}{2},\,0,\,0,\,0,\,0,\,0,\,0,\,0,\,0,\,0,\,-1,\,-1,\,0,\,0,\,0,\,0,\,0,\,0,\,0\right) ,x \rangle \leq 1 \\
-\frac{309485009827534768921207957}{309485009821345068724781056} &\leq \langle \left(0,\,0,\,0,\,0,\,0,\,0,\,0,\,0,\,0,\,1,\,0,\,0,\,0,\,0,\,0,\,0,\,0,\,0,\,0,\,0\right) ,x \rangle \leq \frac{309485009827534768921207957}{309485009821345068724781056} \\
-\frac{9671406563106733594076309}{9671406556917033397649408} &\leq \langle \left(0,\,0,\,0,\,0,\,0,\,0,\,0,\,0,\,0,\,0,\,0,\,1,\,-1,\,0,\,0,\,0,\,0,\,0,\,0,\,0\right) ,x \rangle \leq \frac{4722366500946864110121}{4722366482869645213696} \\
-\frac{40066579960122721}{36028797018963968} &\leq \langle \left(0,\,0,\,0,\,1,\,0,\,0,\,0,\,0,\,0,\,0,\,0,\,0,\,0,\,0,\,0,\,0,\,0,\,0,\,0,\,0\right) ,x \rangle \leq \frac{20040505775302443}{18014398509481984} \\
-\frac{1180593503442084915409}{1180591620717411303424} &\leq \langle \left(0,\,0,\,0,\,0,\,0,\,0,\,0,\,1,\,0,\,0,\,0,\,0,\,0,\,0,\,0,\,0,\,0,\,0,\,0,\,0\right) ,x \rangle \leq \frac{1180593503442084915409}{1180591620717411303424} \\
-\frac{18889473462374272866723}{18889465931478580854784} &\leq \langle \left(0,\,0,\,0,\,0,\,0,\,1,\,0,\,0,\,0,\,0,\,0,\,0,\,0,\,0,\,0,\,0,\,0,\,0,\,0,\,0\right) ,x \rangle \leq \frac{18889473462375840754951}{18889465931478580854784} \\
-\frac{2417851643085731713983019}{2417851639229258349412352} &\leq \langle \left(0,\,0,\,0,\,0,\,0,\,0,\,0,\,0,\,1,\,0,\,0,\,0,\,0,\,0,\,0,\,0,\,0,\,0,\,0,\,0\right) ,x \rangle \leq \frac{2417851643085731713983019}{2417851639229258349412352} \\
-\frac{2361185035008561037417}{2361183241434822606848} &\leq \langle \left(0,\,0,\,0,\,0,\,0,\,0,\,0,\,0,\,0,\,0,\,0,\,0,\,0,\,0,\,1,\,0,\,0,\,0,\,0,\,0\right) ,x \rangle \leq \frac{2361185035008561037417}{2361183241434822606848} \\
-\frac{582225365104347401}{576460752303423488} &\leq \langle \left(1,\,0,\,0,\,0,\,0,\,0,\,0,\,0,\,0,\,0,\,0,\,0,\,100,\,0,\,0,\,0,\,0,\,0,\,0,\,0\right) ,x \rangle \leq \frac{582225507504995595}{576460752303423488} \\
-\frac{2392514266058385}{2251799813685248} &\leq \langle \left(0,\,0,\,0,\,0,\,1,\,0,\,0,\,0,\,0,\,0,\,0,\,0,\,0,\,0,\,0,\,0,\,0,\,0,\,0,\,0\right) ,x \rangle \leq \frac{76542601517243887}{72057594037927936} \\
\end{align*}}
\subsection{Counterexample based on inequality $25$}
If $H$ is the polyhedron obtained from $N$ by removing inequality $I_{25}$, then $H$ is a bounded Hirsch polytope and $N = H \cap C$, where $C$ is the combinatorial $20$-cube defined by
{\footnotesize\begin{align*}
-\frac{5862916589385895}{1125899906842624} &\leq \langle \left(100,\,17,\,-18,\,0,\,0,\,0,\,0,\,0,\,0,\,0,\,0,\,0,\,-1,\,0,\,0,\,0,\,0,\,0,\,0,\,0\right) ,x \rangle \leq 1 \\
-\frac{1180593414291149733993}{1180591620717411303424} &\leq \langle \left(0,\,0,\,0,\,0,\,0,\,0,\,0,\,0,\,0,\,0,\,0,\,0,\,0,\,0,\,0,\,1,\,0,\,0,\,0,\,0\right) ,x \rangle \leq \frac{1180593414291149733993}{1180591620717411303424} \\
-\frac{2361185124158778945083}{2361183241434822606848} &\leq \langle \left(0,\,0,\,0,\,0,\,0,\,0,\,1,\,0,\,0,\,0,\,0,\,0,\,0,\,0,\,0,\,0,\,0,\,0,\,0,\,0\right) ,x \rangle \leq \frac{2361185124159496218833}{2361183241434822606848} \\
-\frac{151699560131112911}{144115188075855872} &\leq \langle \left(0,\,0,\,1,\,0,\,0,\,0,\,0,\,0,\,0,\,0,\,0,\,0,\,-18,\,0,\,0,\,0,\,0,\,0,\,0,\,0\right) ,x \rangle \leq \frac{80781273708265951}{72057594037927936} \\
-\frac{4722372222305608191517}{4722366482869645213696} &\leq \langle \left(0,\,0,\,0,\,0,\,0,\,0,\,0,\,0,\,0,\,0,\,0,\,0,\,0,\,0,\,0,\,0,\,0,\,0,\,1,\,0\right) ,x \rangle \leq \frac{4722372222305608191517}{4722366482869645213696} \\
-\frac{18889470523027351237041}{18889465931478580854784} &\leq \langle \left(0,\,0,\,0,\,0,\,0,\,0,\,0,\,0,\,0,\,0,\,0,\,0,\,0,\,0,\,0,\,0,\,0,\,0,\,0,\,1\right) ,x \rangle \leq \frac{18889470523027351237041}{18889465931478580854784} \\
-\frac{590297603932444082281}{590295810358705651712} &\leq \langle \left(0,\,0,\,0,\,0,\,0,\,0,\,0,\,0,\,0,\,0,\,0,\,0,\,0,\,0,\,0,\,0,\,1,\,0,\,0,\,0\right) ,x \rangle \leq \frac{590297603932444082281}{590295810358705651712} \\
-\frac{9444738705175253405213}{9444732965739290427392} &\leq \langle \left(0,\,0,\,0,\,0,\,0,\,0,\,0,\,0,\,0,\,0,\,0,\,0,\,0,\,0,\,0,\,0,\,0,\,1,\,0,\,0\right) ,x \rangle \leq \frac{9444738705175253405213}{9444732965739290427392} \\
-\frac{39614081257140112583132697033}{39614081257132168796771975168} &\leq \langle \left(0,\,0,\,0,\,0,\,0,\,0,\,0,\,0,\,0,\,0,\,1,\,0,\,0,\,0,\,0,\,0,\,0,\,0,\,0,\,0\right) ,x \rangle \leq \frac{9903520314285028145783175131}{9903520314283042199192993792} \\
-\frac{76859732696219365}{72057594037927936} &\leq \langle \left(0,\,1,\,0,\,0,\,0,\,0,\,0,\,0,\,0,\,0,\,0,\,0,\,17,\,0,\,0,\,0,\,0,\,0,\,0,\,0\right) ,x \rangle \leq \frac{77296343773109141}{72057594037927936} \\
-\frac{147579558414782802697}{147573952589676412928} &\leq \langle \left(0,\,0,\,0,\,0,\,0,\,0,\,0,\,0,\,0,\,0,\,0,\,0,\,0,\,1,\,0,\,0,\,0,\,0,\,0,\,0\right) ,x \rangle \leq \frac{147579558414782802697}{147573952589676412928} \\
-\frac{9671406563106733594076309}{9671406556917033397649408} &\leq \langle \left(0,\,0,\,0,\,0,\,0,\,0,\,0,\,0,\,0,\,0,\,0,\,1,\,0,\,0,\,0,\,0,\,0,\,0,\,0,\,0\right) ,x \rangle \leq \frac{9671406563106733594076309}{9671406556917033397649408} \\
-\frac{309485009827534768921207957}{309485009821345068724781056} &\leq \langle \left(0,\,0,\,0,\,0,\,0,\,0,\,0,\,0,\,0,\,1,\,0,\,0,\,0,\,0,\,0,\,0,\,0,\,0,\,0,\,0\right) ,x \rangle \leq \frac{309485009827534768921207957}{309485009821345068724781056} \\
-\frac{40066579960122721}{36028797018963968} &\leq \langle \left(0,\,0,\,0,\,1,\,0,\,0,\,0,\,0,\,0,\,0,\,0,\,0,\,0,\,0,\,0,\,0,\,0,\,0,\,0,\,0\right) ,x \rangle \leq \frac{20040505775302443}{18014398509481984} \\
-\frac{1180593503442084915409}{1180591620717411303424} &\leq \langle \left(0,\,0,\,0,\,0,\,0,\,0,\,0,\,1,\,0,\,0,\,0,\,0,\,0,\,0,\,0,\,0,\,0,\,0,\,0,\,0\right) ,x \rangle \leq \frac{1180593503442084915409}{1180591620717411303424} \\
-\frac{18889473462374272866723}{18889465931478580854784} &\leq \langle \left(0,\,0,\,0,\,0,\,0,\,1,\,0,\,0,\,0,\,0,\,0,\,0,\,0,\,0,\,0,\,0,\,0,\,0,\,0,\,0\right) ,x \rangle \leq \frac{18889473462375840754951}{18889465931478580854784} \\
-\frac{2417851643085731713983019}{2417851639229258349412352} &\leq \langle \left(0,\,0,\,0,\,0,\,0,\,0,\,0,\,0,\,1,\,0,\,0,\,0,\,0,\,0,\,0,\,0,\,0,\,0,\,0,\,0\right) ,x \rangle \leq \frac{2417851643085731713983019}{2417851639229258349412352} \\
-\frac{2361185035008561037417}{2361183241434822606848} &\leq \langle \left(0,\,0,\,0,\,0,\,0,\,0,\,0,\,0,\,0,\,0,\,0,\,0,\,0,\,0,\,1,\,0,\,0,\,0,\,0,\,0\right) ,x \rangle \leq \frac{2361185035008561037417}{2361183241434822606848} \\
-\frac{582225365104347401}{576460752303423488} &\leq \langle \left(1,\,0,\,0,\,0,\,0,\,0,\,0,\,0,\,0,\,0,\,0,\,0,\,100,\,0,\,0,\,0,\,0,\,0,\,0,\,0\right) ,x \rangle \leq \frac{582225507504995595}{576460752303423488} \\
-\frac{2392514266058385}{2251799813685248} &\leq \langle \left(0,\,0,\,0,\,0,\,1,\,0,\,0,\,0,\,0,\,0,\,0,\,0,\,0,\,0,\,0,\,0,\,0,\,0,\,0,\,0\right) ,x \rangle \leq \frac{76542601517243887}{72057594037927936} \\
\end{align*}}
\subsection{Counterexample based on inequality $26$}
If $H$ is the polyhedron obtained from $N$ by removing inequality $I_{26}$, then $H$ is a bounded Hirsch polytope and $N = H \cap C$, where $C$ is the combinatorial $20$-cube defined by
{\footnotesize\begin{align*}
-\frac{9671406563106733594076309}{9671406556917033397649408} &\leq \langle \left(0,\,0,\,0,\,0,\,0,\,0,\,0,\,0,\,0,\,0,\,0,\,1,\,0,\,0,\,0,\,0,\,0,\,0,\,0,\,0\right) ,x \rangle \leq \frac{9671406563106733594076309}{9671406556917033397649408} \\
-\frac{2361185124158778945083}{2361183241434822606848} &\leq \langle \left(0,\,0,\,0,\,0,\,0,\,0,\,1,\,0,\,0,\,0,\,0,\,0,\,0,\,0,\,0,\,0,\,0,\,0,\,0,\,0\right) ,x \rangle \leq \frac{2361185124159496218833}{2361183241434822606848} \\
-\frac{4722372222305608191517}{4722366482869645213696} &\leq \langle \left(0,\,0,\,0,\,0,\,0,\,0,\,0,\,0,\,0,\,0,\,0,\,0,\,0,\,0,\,0,\,0,\,0,\,0,\,1,\,0\right) ,x \rangle \leq \frac{4722372222305608191517}{4722366482869645213696} \\
-\frac{76859732694334399}{72057594037927936} &\leq \langle \left(0,\,1,\,0,\,0,\,0,\,0,\,0,\,0,\,0,\,0,\,0,\,0,\,0,\,0,\,0,\,0,\,0,\,0,\,0,\,0\right) ,x \rangle \leq \frac{77296343773163779}{72057594037927936} \\
-\frac{18889470523027351237041}{18889465931478580854784} &\leq \langle \left(0,\,0,\,0,\,0,\,0,\,0,\,0,\,0,\,0,\,0,\,0,\,0,\,0,\,0,\,0,\,0,\,0,\,0,\,0,\,1\right) ,x \rangle \leq \frac{18889470523027351237041}{18889465931478580854784} \\
-\frac{590297603932444082281}{590295810358705651712} &\leq \langle \left(0,\,0,\,0,\,0,\,0,\,0,\,0,\,0,\,0,\,0,\,0,\,0,\,0,\,0,\,0,\,0,\,1,\,0,\,0,\,0\right) ,x \rangle \leq \frac{590297603932444082281}{590295810358705651712} \\
-\frac{9444738705175253405213}{9444732965739290427392} &\leq \langle \left(0,\,0,\,0,\,0,\,0,\,0,\,0,\,0,\,0,\,0,\,0,\,0,\,0,\,0,\,0,\,0,\,0,\,1,\,0,\,0\right) ,x \rangle \leq \frac{9444738705175253405213}{9444732965739290427392} \\
-\frac{39614081257140112583132697033}{39614081257132168796771975168} &\leq \langle \left(0,\,0,\,0,\,0,\,0,\,0,\,0,\,0,\,0,\,0,\,1,\,0,\,0,\,0,\,0,\,0,\,0,\,0,\,0,\,0\right) ,x \rangle \leq \frac{9903520314285028145783175131}{9903520314283042199192993792} \\
-\frac{2417851645399615732725419}{2417851639229258349412352} &\leq \langle \left(0,\,0,\,0,\,0,\,0,\,0,\,0,\,0,\,0,\,0,\,0,\,0,\,1,\,0,\,0,\,0,\,0,\,0,\,0,\,0\right) ,x \rangle \leq \frac{2417851645399615732725419}{2417851639229258349412352} \\
-\frac{147579558414782802697}{147573952589676412928} &\leq \langle \left(0,\,0,\,0,\,0,\,0,\,0,\,0,\,0,\,0,\,0,\,0,\,0,\,0,\,1,\,0,\,0,\,0,\,0,\,0,\,0\right) ,x \rangle \leq \frac{147579558414782802697}{147573952589676412928} \\
-\frac{2361185035008561037417}{2361183241434822606848} &\leq \langle \left(0,\,0,\,0,\,0,\,0,\,0,\,0,\,0,\,0,\,0,\,0,\,0,\,0,\,0,\,1,\,0,\,0,\,0,\,0,\,0\right) ,x \rangle \leq \frac{2361185035008561037417}{2361183241434822606848} \\
-\frac{18889473462374272866723}{18889465931478580854784} &\leq \langle \left(0,\,0,\,0,\,0,\,0,\,1,\,0,\,0,\,0,\,0,\,0,\,0,\,0,\,0,\,0,\,0,\,0,\,0,\,0,\,0\right) ,x \rangle \leq \frac{18889473462375840754951}{18889465931478580854784} \\
-\frac{309485009827534768921207957}{309485009821345068724781056} &\leq \langle \left(0,\,0,\,0,\,0,\,0,\,0,\,0,\,0,\,0,\,1,\,0,\,0,\,0,\,0,\,0,\,0,\,0,\,0,\,0,\,0\right) ,x \rangle \leq \frac{309485009827534768921207957}{309485009821345068724781056} \\
-\frac{1180593414291149733993}{1180591620717411303424} &\leq \langle \left(0,\,0,\,0,\,0,\,0,\,0,\,0,\,0,\,0,\,0,\,0,\,0,\,0,\,0,\,0,\,1,\,0,\,0,\,0,\,0\right) ,x \rangle \leq \frac{1180593414291149733993}{1180591620717411303424} \\
-\frac{1180593503442084915409}{1180591620717411303424} &\leq \langle \left(0,\,0,\,0,\,0,\,0,\,0,\,0,\,1,\,0,\,0,\,0,\,0,\,0,\,0,\,0,\,0,\,0,\,0,\,0,\,0\right) ,x \rangle \leq \frac{1180593503442084915409}{1180591620717411303424} \\
-\frac{40066579960122721}{36028797018963968} &\leq \langle \left(0,\,0,\,0,\,1,\,0,\,0,\,0,\,0,\,0,\,0,\,0,\,0,\,0,\,0,\,0,\,0,\,0,\,0,\,0,\,0\right) ,x \rangle \leq \frac{20040505775302443}{18014398509481984} \\
-\frac{2417851643085731713983019}{2417851639229258349412352} &\leq \langle \left(0,\,0,\,0,\,0,\,0,\,0,\,0,\,0,\,1,\,0,\,0,\,0,\,0,\,0,\,0,\,0,\,0,\,0,\,0,\,0\right) ,x \rangle \leq \frac{2417851643085731713983019}{2417851639229258349412352} \\
-\frac{8605210623250955}{2251799813685248} &\leq \langle \left(19,\,0,\,50,\,0,\,0,\,0,\,0,\,0,\,0,\,0,\,0,\,0,\,0,\,0,\,0,\,0,\,0,\,0,\,0,\,0\right) ,x \rangle \leq \frac{7727871631784129}{1125899906842624} \\
-\frac{7428758465248451}{1125899906842624} &\leq \langle \left(100,\,0,\,-38,\,0,\,0,\,0,\,0,\,0,\,0,\,0,\,0,\,0,\,0,\,0,\,0,\,0,\,0,\,0,\,0,\,0\right) ,x \rangle \leq 1 \\
-\frac{2392514266058385}{2251799813685248} &\leq \langle \left(0,\,0,\,0,\,0,\,1,\,0,\,0,\,0,\,0,\,0,\,0,\,0,\,0,\,0,\,0,\,0,\,0,\,0,\,0,\,0\right) ,x \rangle \leq \frac{76542601517243887}{72057594037927936} \\
\end{align*}}
\subsection{Counterexample based on inequality $27$}
If $H$ is the polyhedron obtained from $N$ by removing inequality $I_{27}$, then $H$ is a bounded Hirsch polytope and $N = H \cap C$, where $C$ is the combinatorial $20$-cube defined by
{\footnotesize\begin{align*}
-\frac{6351125391016935}{1125899906842624} &\leq \langle \left(100,\,-22,\,-17,\,0,\,0,\,0,\,0,\,0,\,0,\,0,\,0,\,0,\,0,\,0,\,0,\,0,\,0,\,0,\,0,\,0\right) ,x \rangle \leq 1 \\
-\frac{1180593414291149733993}{1180591620717411303424} &\leq \langle \left(0,\,0,\,0,\,0,\,0,\,0,\,0,\,0,\,0,\,0,\,0,\,0,\,0,\,0,\,0,\,1,\,0,\,0,\,0,\,0\right) ,x \rangle \leq \frac{1180593414291149733993}{1180591620717411303424} \\
-\frac{9671406563106733594076309}{9671406556917033397649408} &\leq \langle \left(0,\,0,\,0,\,0,\,0,\,0,\,0,\,0,\,0,\,0,\,0,\,1,\,0,\,0,\,0,\,0,\,0,\,0,\,0,\,0\right) ,x \rangle \leq \frac{9671406563106733594076309}{9671406556917033397649408} \\
-\frac{1180593503442084915409}{1180591620717411303424} &\leq \langle \left(0,\,0,\,0,\,0,\,0,\,0,\,0,\,1,\,0,\,0,\,0,\,0,\,0,\,0,\,0,\,0,\,0,\,0,\,0,\,0\right) ,x \rangle \leq \frac{1180593503442084915409}{1180591620717411303424} \\
-\frac{4722372222305608191517}{4722366482869645213696} &\leq \langle \left(0,\,0,\,0,\,0,\,0,\,0,\,0,\,0,\,0,\,0,\,0,\,0,\,0,\,0,\,0,\,0,\,0,\,0,\,1,\,0\right) ,x \rangle \leq \frac{4722372222305608191517}{4722366482869645213696} \\
-\frac{18889470523027351237041}{18889465931478580854784} &\leq \langle \left(0,\,0,\,0,\,0,\,0,\,0,\,0,\,0,\,0,\,0,\,0,\,0,\,0,\,0,\,0,\,0,\,0,\,0,\,0,\,1\right) ,x \rangle \leq \frac{18889470523027351237041}{18889465931478580854784} \\
-\frac{590297603932444082281}{590295810358705651712} &\leq \langle \left(0,\,0,\,0,\,0,\,0,\,0,\,0,\,0,\,0,\,0,\,0,\,0,\,0,\,0,\,0,\,0,\,1,\,0,\,0,\,0\right) ,x \rangle \leq \frac{590297603932444082281}{590295810358705651712} \\
-\frac{9444738705175253405213}{9444732965739290427392} &\leq \langle \left(0,\,0,\,0,\,0,\,0,\,0,\,0,\,0,\,0,\,0,\,0,\,0,\,0,\,0,\,0,\,0,\,0,\,1,\,0,\,0\right) ,x \rangle \leq \frac{9444738705175253405213}{9444732965739290427392} \\
-\frac{11462966946396265}{9007199254740992} &\leq \langle \left(1,\,3,\,2,\,0,\,0,\,0,\,0,\,0,\,0,\,0,\,0,\,0,\,0,\,0,\,0,\,0,\,0,\,0,\,0,\,0\right) ,x \rangle \leq \frac{6523746167780533}{4503599627370496} \\
-\frac{39614081257140112583132697033}{39614081257132168796771975168} &\leq \langle \left(0,\,0,\,0,\,0,\,0,\,0,\,0,\,0,\,0,\,0,\,1,\,0,\,0,\,0,\,0,\,0,\,0,\,0,\,0,\,0\right) ,x \rangle \leq \frac{9903520314285028145783175131}{9903520314283042199192993792} \\
-\frac{2417851645399615732725419}{2417851639229258349412352} &\leq \langle \left(0,\,0,\,0,\,0,\,0,\,0,\,0,\,0,\,0,\,0,\,0,\,0,\,1,\,0,\,0,\,0,\,0,\,0,\,0,\,0\right) ,x \rangle \leq \frac{2417851645399615732725419}{2417851639229258349412352} \\
-\frac{147579558414782802697}{147573952589676412928} &\leq \langle \left(0,\,0,\,0,\,0,\,0,\,0,\,0,\,0,\,0,\,0,\,0,\,0,\,0,\,1,\,0,\,0,\,0,\,0,\,0,\,0\right) ,x \rangle \leq \frac{147579558414782802697}{147573952589676412928} \\
-\frac{18889473462374272866723}{18889465931478580854784} &\leq \langle \left(0,\,0,\,0,\,0,\,0,\,1,\,0,\,0,\,0,\,0,\,0,\,0,\,0,\,0,\,0,\,0,\,0,\,0,\,0,\,0\right) ,x \rangle \leq \frac{18889473462375840754951}{18889465931478580854784} \\
-\frac{309485009827534768921207957}{309485009821345068724781056} &\leq \langle \left(0,\,0,\,0,\,0,\,0,\,0,\,0,\,0,\,0,\,1,\,0,\,0,\,0,\,0,\,0,\,0,\,0,\,0,\,0,\,0\right) ,x \rangle \leq \frac{309485009827534768921207957}{309485009821345068724781056} \\
-\frac{40066579960122721}{36028797018963968} &\leq \langle \left(0,\,0,\,0,\,1,\,0,\,0,\,0,\,0,\,0,\,0,\,0,\,0,\,0,\,0,\,0,\,0,\,0,\,0,\,0,\,0\right) ,x \rangle \leq \frac{20040505775302443}{18014398509481984} \\
-\frac{5259431040373345}{1125899906842624} &\leq \langle \left(0,\,17,\,-22,\,0,\,0,\,0,\,0,\,0,\,0,\,0,\,0,\,0,\,0,\,0,\,0,\,0,\,0,\,0,\,0,\,0\right) ,x \rangle \leq \frac{1807107238614929}{562949953421312} \\
-\frac{2361185124158778945083}{2361183241434822606848} &\leq \langle \left(0,\,0,\,0,\,0,\,0,\,0,\,1,\,0,\,0,\,0,\,0,\,0,\,0,\,0,\,0,\,0,\,0,\,0,\,0,\,0\right) ,x \rangle \leq \frac{2361185124159496218833}{2361183241434822606848} \\
-\frac{2417851643085731713983019}{2417851639229258349412352} &\leq \langle \left(0,\,0,\,0,\,0,\,0,\,0,\,0,\,0,\,1,\,0,\,0,\,0,\,0,\,0,\,0,\,0,\,0,\,0,\,0,\,0\right) ,x \rangle \leq \frac{2417851643085731713983019}{2417851639229258349412352} \\
-\frac{2361185035008561037417}{2361183241434822606848} &\leq \langle \left(0,\,0,\,0,\,0,\,0,\,0,\,0,\,0,\,0,\,0,\,0,\,0,\,0,\,0,\,1,\,0,\,0,\,0,\,0,\,0\right) ,x \rangle \leq \frac{2361185035008561037417}{2361183241434822606848} \\
-\frac{2392514266058385}{2251799813685248} &\leq \langle \left(0,\,0,\,0,\,0,\,1,\,0,\,0,\,0,\,0,\,0,\,0,\,0,\,0,\,0,\,0,\,0,\,0,\,0,\,0,\,0\right) ,x \rangle \leq \frac{76542601517243887}{72057594037927936} \\
\end{align*}}
\subsection{Counterexample based on inequality $28$}
If $H$ is the polyhedron obtained from $N$ by removing inequality $I_{28}$, then $H$ is a bounded Hirsch polytope and $N = H \cap C$, where $C$ is the combinatorial $20$-cube defined by
{\footnotesize\begin{align*}
-\frac{1180593414291149733993}{1180591620717411303424} &\leq \langle \left(0,\,0,\,0,\,0,\,0,\,0,\,0,\,0,\,0,\,0,\,0,\,0,\,0,\,0,\,0,\,1,\,0,\,0,\,0,\,0\right) ,x \rangle \leq \frac{1180593414291149733993}{1180591620717411303424} \\
-\frac{9671406563106733594076309}{9671406556917033397649408} &\leq \langle \left(0,\,0,\,0,\,0,\,0,\,0,\,0,\,0,\,0,\,0,\,0,\,1,\,0,\,0,\,0,\,0,\,0,\,0,\,0,\,0\right) ,x \rangle \leq \frac{9671406563106733594076309}{9671406556917033397649408} \\
-\frac{1180593503442084915409}{1180591620717411303424} &\leq \langle \left(0,\,0,\,0,\,0,\,0,\,0,\,0,\,1,\,0,\,0,\,0,\,0,\,0,\,0,\,0,\,0,\,0,\,0,\,0,\,0\right) ,x \rangle \leq \frac{1180593503442084915409}{1180591620717411303424} \\
-\frac{4722372222305608191517}{4722366482869645213696} &\leq \langle \left(0,\,0,\,0,\,0,\,0,\,0,\,0,\,0,\,0,\,0,\,0,\,0,\,0,\,0,\,0,\,0,\,0,\,0,\,1,\,0\right) ,x \rangle \leq \frac{4722372222305608191517}{4722366482869645213696} \\
-\frac{18889470523027351237041}{18889465931478580854784} &\leq \langle \left(0,\,0,\,0,\,0,\,0,\,0,\,0,\,0,\,0,\,0,\,0,\,0,\,0,\,0,\,0,\,0,\,0,\,0,\,0,\,1\right) ,x \rangle \leq \frac{18889470523027351237041}{18889465931478580854784} \\
-\frac{2417851645399615732725419}{2417851639229258349412352} &\leq \langle \left(0,\,0,\,0,\,0,\,0,\,0,\,0,\,0,\,0,\,0,\,0,\,0,\,1,\,0,\,0,\,0,\,0,\,0,\,0,\,0\right) ,x \rangle \leq \frac{2417851645399615732725419}{2417851639229258349412352} \\
-\frac{590297603932444082281}{590295810358705651712} &\leq \langle \left(0,\,0,\,0,\,0,\,0,\,0,\,0,\,0,\,0,\,0,\,0,\,0,\,0,\,0,\,0,\,0,\,1,\,0,\,0,\,0\right) ,x \rangle \leq \frac{590297603932444082281}{590295810358705651712} \\
-\frac{9444738705175253405213}{9444732965739290427392} &\leq \langle \left(0,\,0,\,0,\,0,\,0,\,0,\,0,\,0,\,0,\,0,\,0,\,0,\,0,\,0,\,0,\,0,\,0,\,1,\,0,\,0\right) ,x \rangle \leq \frac{9444738705175253405213}{9444732965739290427392} \\
-\frac{42309792705210917}{36028797018963968} &\leq \langle \left(0,\,0,\,0,\,1,\,1,\,0,\,0,\,0,\,0,\,0,\,0,\,0,\,0,\,0,\,0,\,0,\,0,\,0,\,0,\,0\right) ,x \rangle \leq \frac{20324969558447131}{18014398509481984} \\
-\frac{39614081257140112583132697033}{39614081257132168796771975168} &\leq \langle \left(0,\,0,\,0,\,0,\,0,\,0,\,0,\,0,\,0,\,0,\,1,\,0,\,0,\,0,\,0,\,0,\,0,\,0,\,0,\,0\right) ,x \rangle \leq \frac{9903520314285028145783175131}{9903520314283042199192993792} \\
-\frac{9038463700167585}{281474976710656} &\leq \langle \left(1,\,0,\,0,\,0,\,-500,\,0,\,0,\,0,\,0,\,0,\,0,\,0,\,0,\,0,\,0,\,0,\,0,\,0,\,0,\,0\right) ,x \rangle \leq \frac{9078942084466643}{281474976710656} \\
-\frac{147579558414782802697}{147573952589676412928} &\leq \langle \left(0,\,0,\,0,\,0,\,0,\,0,\,0,\,0,\,0,\,0,\,0,\,0,\,0,\,1,\,0,\,0,\,0,\,0,\,0,\,0\right) ,x \rangle \leq \frac{147579558414782802697}{147573952589676412928} \\
-\frac{18889473462374272866723}{18889465931478580854784} &\leq \langle \left(0,\,0,\,0,\,0,\,0,\,1,\,0,\,0,\,0,\,0,\,0,\,0,\,0,\,0,\,0,\,0,\,0,\,0,\,0,\,0\right) ,x \rangle \leq \frac{18889473462375840754951}{18889465931478580854784} \\
-\frac{309485009827534768921207957}{309485009821345068724781056} &\leq \langle \left(0,\,0,\,0,\,0,\,0,\,0,\,0,\,0,\,0,\,1,\,0,\,0,\,0,\,0,\,0,\,0,\,0,\,0,\,0,\,0\right) ,x \rangle \leq \frac{309485009827534768921207957}{309485009821345068724781056} \\
-\frac{151699951751958945}{144115188075855872} &\leq \langle \left(0,\,0,\,1,\,0,\,0,\,0,\,0,\,0,\,0,\,0,\,0,\,0,\,0,\,0,\,0,\,0,\,0,\,0,\,0,\,0\right) ,x \rangle \leq \frac{80781273703537389}{72057594037927936} \\
-\frac{4629497230848269}{1125899906842624} &\leq \langle \left(0,\,1,\,0,\,0,\,-50,\,0,\,0,\,0,\,0,\,0,\,0,\,0,\,0,\,0,\,0,\,0,\,0,\,0,\,0,\,0\right) ,x \rangle \leq \frac{9287544710655669}{2251799813685248} \\
-\frac{2417851643085731713983019}{2417851639229258349412352} &\leq \langle \left(0,\,0,\,0,\,0,\,0,\,0,\,0,\,0,\,1,\,0,\,0,\,0,\,0,\,0,\,0,\,0,\,0,\,0,\,0,\,0\right) ,x \rangle \leq \frac{2417851643085731713983019}{2417851639229258349412352} \\
-\frac{2361185035008561037417}{2361183241434822606848} &\leq \langle \left(0,\,0,\,0,\,0,\,0,\,0,\,0,\,0,\,0,\,0,\,0,\,0,\,0,\,0,\,1,\,0,\,0,\,0,\,0,\,0\right) ,x \rangle \leq \frac{2361185035008561037417}{2361183241434822606848} \\
-\frac{12005402462150447}{4503599627370496} &\leq \langle \left(100,\,10,\,0,\,-\frac{1}{5},\,\frac{1}{5},\,0,\,0,\,0,\,0,\,0,\,0,\,0,\,0,\,0,\,0,\,0,\,0,\,0,\,0,\,0\right) ,x \rangle \leq 1 \\
-\frac{2361185124158778945083}{2361183241434822606848} &\leq \langle \left(0,\,0,\,0,\,0,\,0,\,0,\,1,\,0,\,0,\,0,\,0,\,0,\,0,\,0,\,0,\,0,\,0,\,0,\,0,\,0\right) ,x \rangle \leq \frac{2361185124159496218833}{2361183241434822606848} \\
\end{align*}}
\subsection{Counterexample based on inequality $29$}
If $H$ is the polyhedron obtained from $N$ by removing inequality $I_{29}$, then $H$ is a bounded Hirsch polytope and $N = H \cap C$, where $C$ is the combinatorial $20$-cube defined by
{\footnotesize\begin{align*}
-\frac{2361185124158778945083}{2361183241434822606848} &\leq \langle \left(0,\,0,\,0,\,0,\,0,\,0,\,1,\,0,\,0,\,0,\,0,\,0,\,0,\,0,\,0,\,0,\,0,\,0,\,0,\,0\right) ,x \rangle \leq \frac{2361185124159496218833}{2361183241434822606848} \\
-\frac{20033293414868413}{18014398509481984} &\leq \langle \left(0,\,0,\,0,\,1,\,0,\,0,\,0,\,\frac{3}{25},\,0,\,0,\,0,\,0,\,0,\,0,\,0,\,0,\,0,\,0,\,0,\,0\right) ,x \rangle \leq \frac{80162036889545207}{72057594037927936} \\
-\frac{4722372222305608191517}{4722366482869645213696} &\leq \langle \left(0,\,0,\,0,\,0,\,0,\,0,\,0,\,0,\,0,\,0,\,0,\,0,\,0,\,0,\,0,\,0,\,0,\,0,\,1,\,0\right) ,x \rangle \leq \frac{4722372222305608191517}{4722366482869645213696} \\
-\frac{38280239746676815}{36028797018963968} &\leq \langle \left(0,\,0,\,0,\,0,\,1,\,0,\,0,\,\frac{1}{5},\,0,\,0,\,0,\,0,\,0,\,0,\,0,\,0,\,0,\,0,\,0,\,0\right) ,x \rangle \leq \frac{9567825198625803}{9007199254740992} \\
-\frac{590297603932444082281}{590295810358705651712} &\leq \langle \left(0,\,0,\,0,\,0,\,0,\,0,\,0,\,0,\,0,\,0,\,0,\,0,\,0,\,0,\,0,\,0,\,1,\,0,\,0,\,0\right) ,x \rangle \leq \frac{590297603932444082281}{590295810358705651712} \\
-\frac{18889470523027351237041}{18889465931478580854784} &\leq \langle \left(0,\,0,\,0,\,0,\,0,\,0,\,0,\,0,\,0,\,0,\,0,\,0,\,0,\,0,\,0,\,0,\,0,\,0,\,0,\,1\right) ,x \rangle \leq \frac{18889470523027351237041}{18889465931478580854784} \\
-\frac{291112682757822535}{288230376151711744} &\leq \langle \left(1,\,0,\,0,\,0,\,0,\,0,\,0,\,100,\,0,\,0,\,0,\,0,\,0,\,0,\,0,\,0,\,0,\,0,\,0,\,0\right) ,x \rangle \leq \frac{145578850792339161}{144115188075855872} \\
-\frac{9444738705175253405213}{9444732965739290427392} &\leq \langle \left(0,\,0,\,0,\,0,\,0,\,0,\,0,\,0,\,0,\,0,\,0,\,0,\,0,\,0,\,0,\,0,\,0,\,1,\,0,\,0\right) ,x \rangle \leq \frac{9444738705175253405213}{9444732965739290427392} \\
-\frac{39614081257140112583132697033}{39614081257132168796771975168} &\leq \langle \left(0,\,0,\,0,\,0,\,0,\,0,\,0,\,0,\,0,\,0,\,1,\,0,\,0,\,0,\,0,\,0,\,0,\,0,\,0,\,0\right) ,x \rangle \leq \frac{9903520314285028145783175131}{9903520314283042199192993792} \\
-\frac{2417851645399615732725419}{2417851639229258349412352} &\leq \langle \left(0,\,0,\,0,\,0,\,0,\,0,\,0,\,0,\,0,\,0,\,0,\,0,\,1,\,0,\,0,\,0,\,0,\,0,\,0,\,0\right) ,x \rangle \leq \frac{2417851645399615732725419}{2417851639229258349412352} \\
-\frac{147579558414782802697}{147573952589676412928} &\leq \langle \left(0,\,0,\,0,\,0,\,0,\,0,\,0,\,0,\,0,\,0,\,0,\,0,\,0,\,1,\,0,\,0,\,0,\,0,\,0,\,0\right) ,x \rangle \leq \frac{147579558414782802697}{147573952589676412928} \\
-\frac{76859733766363671}{72057594037927936} &\leq \langle \left(0,\,1,\,0,\,0,\,0,\,0,\,0,\,-\frac{2999}{100},\,0,\,0,\,0,\,0,\,0,\,0,\,0,\,0,\,0,\,0,\,0,\,0\right) ,x \rangle \leq \frac{77296299606963433}{72057594037927936} \\
-\frac{4705868905881905}{1125899906842624} &\leq \langle \left(100,\,-\frac{2999}{100},\,0,\,\frac{3}{25},\,\frac{1}{5},\,0,\,0,\,-1,\,0,\,0,\,0,\,0,\,0,\,0,\,0,\,0,\,0,\,0,\,0,\,0\right) ,x \rangle \leq 1 \\
-\frac{9671406563106733594076309}{9671406556917033397649408} &\leq \langle \left(0,\,0,\,0,\,0,\,0,\,0,\,0,\,0,\,0,\,0,\,0,\,1,\,0,\,0,\,0,\,0,\,0,\,0,\,0,\,0\right) ,x \rangle \leq \frac{9671406563106733594076309}{9671406556917033397649408} \\
-\frac{309485009827534768921207957}{309485009821345068724781056} &\leq \langle \left(0,\,0,\,0,\,0,\,0,\,0,\,0,\,0,\,0,\,1,\,0,\,0,\,0,\,0,\,0,\,0,\,0,\,0,\,0,\,0\right) ,x \rangle \leq \frac{309485009827534768921207957}{309485009821345068724781056} \\
-\frac{1180593414291149733993}{1180591620717411303424} &\leq \langle \left(0,\,0,\,0,\,0,\,0,\,0,\,0,\,0,\,0,\,0,\,0,\,0,\,0,\,0,\,0,\,1,\,0,\,0,\,0,\,0\right) ,x \rangle \leq \frac{1180593414291149733993}{1180591620717411303424} \\
-\frac{151699951751958945}{144115188075855872} &\leq \langle \left(0,\,0,\,1,\,0,\,0,\,0,\,0,\,0,\,0,\,0,\,0,\,0,\,0,\,0,\,0,\,0,\,0,\,0,\,0,\,0\right) ,x \rangle \leq \frac{80781273703537389}{72057594037927936} \\
-\frac{2417851643085731713983019}{2417851639229258349412352} &\leq \langle \left(0,\,0,\,0,\,0,\,0,\,0,\,0,\,0,\,1,\,0,\,0,\,0,\,0,\,0,\,0,\,0,\,0,\,0,\,0,\,0\right) ,x \rangle \leq \frac{2417851643085731713983019}{2417851639229258349412352} \\
-\frac{2361185035008561037417}{2361183241434822606848} &\leq \langle \left(0,\,0,\,0,\,0,\,0,\,0,\,0,\,0,\,0,\,0,\,0,\,0,\,0,\,0,\,1,\,0,\,0,\,0,\,0,\,0\right) ,x \rangle \leq \frac{2361185035008561037417}{2361183241434822606848} \\
-\frac{18889473462374272866723}{18889465931478580854784} &\leq \langle \left(0,\,0,\,0,\,0,\,0,\,1,\,0,\,0,\,0,\,0,\,0,\,0,\,0,\,0,\,0,\,0,\,0,\,0,\,0,\,0\right) ,x \rangle \leq \frac{18889473462375840754951}{18889465931478580854784} \\
\end{align*}}
\subsection{Counterexample based on inequality $30$}
If $H$ is the polyhedron obtained from $N$ by removing inequality $I_{30}$, then $H$ is a bounded Hirsch polytope and $N = H \cap C$, where $C$ is the combinatorial $20$-cube defined by
{\footnotesize\begin{align*}
-\frac{1180593414291149733993}{1180591620717411303424} &\leq \langle \left(0,\,0,\,0,\,0,\,0,\,0,\,0,\,0,\,0,\,0,\,0,\,0,\,0,\,0,\,0,\,1,\,0,\,0,\,0,\,0\right) ,x \rangle \leq \frac{1180593414291149733993}{1180591620717411303424} \\
-\frac{2361185035008561037417}{2361183241434822606848} &\leq \langle \left(0,\,0,\,0,\,0,\,0,\,0,\,0,\,0,\,0,\,0,\,0,\,0,\,0,\,0,\,1,\,0,\,0,\,0,\,0,\,0\right) ,x \rangle \leq \frac{2361185035008561037417}{2361183241434822606848} \\
-\frac{1180593503442084915409}{1180591620717411303424} &\leq \langle \left(0,\,0,\,0,\,0,\,0,\,0,\,0,\,1,\,0,\,0,\,0,\,0,\,0,\,0,\,0,\,0,\,0,\,0,\,0,\,0\right) ,x \rangle \leq \frac{1180593503442084915409}{1180591620717411303424} \\
-\frac{4722372222305608191517}{4722366482869645213696} &\leq \langle \left(0,\,0,\,0,\,0,\,0,\,0,\,0,\,0,\,0,\,0,\,0,\,0,\,0,\,0,\,0,\,0,\,0,\,0,\,1,\,0\right) ,x \rangle \leq \frac{4722372222305608191517}{4722366482869645213696} \\
-\frac{4803733326909015}{4503599627370496} &\leq \langle \left(0,\,1,\,0,\,0,\,0,\,0,\,-\frac{299999}{10000},\,0,\,0,\,0,\,0,\,0,\,0,\,0,\,0,\,0,\,0,\,0,\,0,\,0\right) ,x \rangle \leq \frac{77296299606963433}{72057594037927936} \\
-\frac{18889470523027351237041}{18889465931478580854784} &\leq \langle \left(0,\,0,\,0,\,0,\,0,\,0,\,0,\,0,\,0,\,0,\,0,\,0,\,0,\,0,\,0,\,0,\,0,\,0,\,0,\,1\right) ,x \rangle \leq \frac{18889470523027351237041}{18889465931478580854784} \\
-\frac{590297603932444082281}{590295810358705651712} &\leq \langle \left(0,\,0,\,0,\,0,\,0,\,0,\,0,\,0,\,0,\,0,\,0,\,0,\,0,\,0,\,0,\,0,\,1,\,0,\,0,\,0\right) ,x \rangle \leq \frac{590297603932444082281}{590295810358705651712} \\
-\frac{9444738705175253405213}{9444732965739290427392} &\leq \langle \left(0,\,0,\,0,\,0,\,0,\,0,\,0,\,0,\,0,\,0,\,0,\,0,\,0,\,0,\,0,\,0,\,0,\,1,\,0,\,0\right) ,x \rangle \leq \frac{9444738705175253405213}{9444732965739290427392} \\
-\frac{39614081257140112583132697033}{39614081257132168796771975168} &\leq \langle \left(0,\,0,\,0,\,0,\,0,\,0,\,0,\,0,\,0,\,0,\,1,\,0,\,0,\,0,\,0,\,0,\,0,\,0,\,0,\,0\right) ,x \rangle \leq \frac{9903520314285028145783175131}{9903520314283042199192993792} \\
-\frac{2417851645399615732725419}{2417851639229258349412352} &\leq \langle \left(0,\,0,\,0,\,0,\,0,\,0,\,0,\,0,\,0,\,0,\,0,\,0,\,1,\,0,\,0,\,0,\,0,\,0,\,0,\,0\right) ,x \rangle \leq \frac{2417851645399615732725419}{2417851639229258349412352} \\
-\frac{147579558414782802697}{147573952589676412928} &\leq \langle \left(0,\,0,\,0,\,0,\,0,\,0,\,0,\,0,\,0,\,0,\,0,\,0,\,0,\,1,\,0,\,0,\,0,\,0,\,0,\,0\right) ,x \rangle \leq \frac{147579558414782802697}{147573952589676412928} \\
-\frac{9671406563106733594076309}{9671406556917033397649408} &\leq \langle \left(0,\,0,\,0,\,0,\,0,\,0,\,0,\,0,\,0,\,0,\,0,\,1,\,0,\,0,\,0,\,0,\,0,\,0,\,0,\,0\right) ,x \rangle \leq \frac{9671406563106733594076309}{9671406556917033397649408} \\
-\frac{309485009827534768921207957}{309485009821345068724781056} &\leq \langle \left(0,\,0,\,0,\,0,\,0,\,0,\,0,\,0,\,0,\,1,\,0,\,0,\,0,\,0,\,0,\,0,\,0,\,0,\,0,\,0\right) ,x \rangle \leq \frac{309485009827534768921207957}{309485009821345068724781056} \\
-\frac{4706681271745485}{1125899906842624} &\leq \langle \left(100,\,-\frac{299999}{10000},\,0,\,0,\,-\frac{1}{100},\,0,\,-1,\,0,\,0,\,0,\,0,\,0,\,0,\,0,\,0,\,0,\,0,\,0,\,0,\,0\right) ,x \rangle \leq 1 \\
-\frac{151699951751958945}{144115188075855872} &\leq \langle \left(0,\,0,\,1,\,0,\,0,\,0,\,0,\,0,\,0,\,0,\,0,\,0,\,0,\,0,\,0,\,0,\,0,\,0,\,0,\,0\right) ,x \rangle \leq \frac{80781273703537389}{72057594037927936} \\
-\frac{40066579960122721}{36028797018963968} &\leq \langle \left(0,\,0,\,0,\,1,\,0,\,0,\,0,\,0,\,0,\,0,\,0,\,0,\,0,\,0,\,0,\,0,\,0,\,0,\,0,\,0\right) ,x \rangle \leq \frac{20040505775302443}{18014398509481984} \\
-\frac{2417851643085731713983019}{2417851639229258349412352} &\leq \langle \left(0,\,0,\,0,\,0,\,0,\,0,\,0,\,0,\,1,\,0,\,0,\,0,\,0,\,0,\,0,\,0,\,0,\,0,\,0,\,0\right) ,x \rangle \leq \frac{2417851643085731713983019}{2417851639229258349412352} \\
-\frac{38280228544216787}{36028797018963968} &\leq \langle \left(0,\,0,\,0,\,0,\,1,\,0,\,-\frac{1}{100},\,0,\,0,\,0,\,0,\,0,\,0,\,0,\,0,\,0,\,0,\,0,\,0,\,0\right) ,x \rangle \leq \frac{153085203030899647}{144115188075855872} \\
-\frac{1137158916616491}{1125899906842624} &\leq \langle \left(1,\,0,\,0,\,0,\,0,\,0,\,100,\,0,\,0,\,0,\,0,\,0,\,0,\,0,\,0,\,0,\,0,\,0,\,0,\,0\right) ,x \rangle \leq \frac{291134718970022493}{288230376151711744} \\
-\frac{18889473462374272866723}{18889465931478580854784} &\leq \langle \left(0,\,0,\,0,\,0,\,0,\,1,\,0,\,0,\,0,\,0,\,0,\,0,\,0,\,0,\,0,\,0,\,0,\,0,\,0,\,0\right) ,x \rangle \leq \frac{18889473462375840754951}{18889465931478580854784} \\
\end{align*}}
\subsection{Counterexample based on inequality $31$}
If $H$ is the polyhedron obtained from $N$ by removing inequality $I_{31}$, then $H$ is a bounded Hirsch polytope and $N = H \cap C$, where $C$ is the combinatorial $20$-cube defined by
{\footnotesize\begin{align*}
-\frac{18889470523027351237041}{18889465931478580854784} &\leq \langle \left(0,\,0,\,0,\,0,\,0,\,0,\,0,\,0,\,0,\,0,\,0,\,0,\,0,\,0,\,0,\,0,\,0,\,0,\,0,\,1\right) ,x \rangle \leq \frac{18889470523027351237041}{18889465931478580854784} \\
-\frac{8621339511414579}{2251799813685248} &\leq \langle \left(100,\,\frac{549}{20},\,0,\,-\frac{1}{5000},\,-\frac{1}{800},\,-1,\,0,\,0,\,0,\,0,\,0,\,0,\,0,\,0,\,0,\,0,\,0,\,0,\,0,\,0\right) ,x \rangle \leq 1 \\
-\frac{9671406563106733594076309}{9671406556917033397649408} &\leq \langle \left(0,\,0,\,0,\,0,\,0,\,0,\,0,\,0,\,0,\,0,\,0,\,1,\,0,\,0,\,0,\,0,\,0,\,0,\,0,\,0\right) ,x \rangle \leq \frac{9671406563106733594076309}{9671406556917033397649408} \\
-\frac{1180593503442084915409}{1180591620717411303424} &\leq \langle \left(0,\,0,\,0,\,0,\,0,\,0,\,0,\,1,\,0,\,0,\,0,\,0,\,0,\,0,\,0,\,0,\,0,\,0,\,0,\,0\right) ,x \rangle \leq \frac{1180593503442084915409}{1180591620717411303424} \\
-\frac{2417851643085731713983019}{2417851639229258349412352} &\leq \langle \left(0,\,0,\,0,\,0,\,0,\,0,\,0,\,0,\,1,\,0,\,0,\,0,\,0,\,0,\,0,\,0,\,0,\,0,\,0,\,0\right) ,x \rangle \leq \frac{2417851643085731713983019}{2417851639229258349412352} \\
-\frac{590297603932444082281}{590295810358705651712} &\leq \langle \left(0,\,0,\,0,\,0,\,0,\,0,\,0,\,0,\,0,\,0,\,0,\,0,\,0,\,0,\,0,\,0,\,1,\,0,\,0,\,0\right) ,x \rangle \leq \frac{590297603932444082281}{590295810358705651712} \\
-\frac{2417851645399615732725419}{2417851639229258349412352} &\leq \langle \left(0,\,0,\,0,\,0,\,0,\,0,\,0,\,0,\,0,\,0,\,0,\,0,\,1,\,0,\,0,\,0,\,0,\,0,\,0,\,0\right) ,x \rangle \leq \frac{2417851645399615732725419}{2417851639229258349412352} \\
-\frac{582225365202441111}{576460752303423488} &\leq \langle \left(1,\,0,\,0,\,0,\,0,\,100,\,0,\,0,\,0,\,0,\,0,\,0,\,0,\,0,\,0,\,0,\,0,\,0,\,0,\,0\right) ,x \rangle \leq \frac{582246455329420841}{576460752303423488} \\
-\frac{9444738705175253405213}{9444732965739290427392} &\leq \langle \left(0,\,0,\,0,\,0,\,0,\,0,\,0,\,0,\,0,\,0,\,0,\,0,\,0,\,0,\,0,\,0,\,0,\,1,\,0,\,0\right) ,x \rangle \leq \frac{9444738705175253405213}{9444732965739290427392} \\
-\frac{39614081257140112583132697033}{39614081257132168796771975168} &\leq \langle \left(0,\,0,\,0,\,0,\,0,\,0,\,0,\,0,\,0,\,0,\,1,\,0,\,0,\,0,\,0,\,0,\,0,\,0,\,0,\,0\right) ,x \rangle \leq \frac{9903520314285028145783175131}{9903520314283042199192993792} \\
-\frac{10016644990746265}{9007199254740992} &\leq \langle \left(0,\,0,\,0,\,1,\,0,\,-\frac{1}{5000},\,0,\,0,\,0,\,0,\,0,\,0,\,0,\,0,\,0,\,0,\,0,\,0,\,0,\,0\right) ,x \rangle \leq \frac{10020252886933079}{9007199254740992} \\
-\frac{147579558414782802697}{147573952589676412928} &\leq \langle \left(0,\,0,\,0,\,0,\,0,\,0,\,0,\,0,\,0,\,0,\,0,\,0,\,0,\,1,\,0,\,0,\,0,\,0,\,0,\,0\right) ,x \rangle \leq \frac{147579558414782802697}{147573952589676412928} \\
-\frac{76859732938478303}{72057594037927936} &\leq \langle \left(0,\,1,\,0,\,0,\,0,\,\frac{549}{20},\,0,\,0,\,0,\,0,\,0,\,0,\,0,\,0,\,0,\,0,\,0,\,0,\,0,\,0\right) ,x \rangle \leq \frac{38648171883146663}{36028797018963968} \\
-\frac{309485009827534768921207957}{309485009821345068724781056} &\leq \langle \left(0,\,0,\,0,\,0,\,0,\,0,\,0,\,0,\,0,\,1,\,0,\,0,\,0,\,0,\,0,\,0,\,0,\,0,\,0,\,0\right) ,x \rangle \leq \frac{309485009827534768921207957}{309485009821345068724781056} \\
-\frac{1180593414291149733993}{1180591620717411303424} &\leq \langle \left(0,\,0,\,0,\,0,\,0,\,0,\,0,\,0,\,0,\,0,\,0,\,0,\,0,\,0,\,0,\,1,\,0,\,0,\,0,\,0\right) ,x \rangle \leq \frac{1180593414291149733993}{1180591620717411303424} \\
-\frac{151699951751958945}{144115188075855872} &\leq \langle \left(0,\,0,\,1,\,0,\,0,\,0,\,0,\,0,\,0,\,0,\,0,\,0,\,0,\,0,\,0,\,0,\,0,\,0,\,0,\,0\right) ,x \rangle \leq \frac{80781273703537389}{72057594037927936} \\
-\frac{153120913099557297}{144115188075855872} &\leq \langle \left(0,\,0,\,0,\,0,\,1,\,-\frac{1}{800},\,0,\,0,\,0,\,0,\,0,\,0,\,0,\,0,\,0,\,0,\,0,\,0,\,0,\,0\right) ,x \rangle \leq \frac{38271300758565879}{36028797018963968} \\
-\frac{2361185035008561037417}{2361183241434822606848} &\leq \langle \left(0,\,0,\,0,\,0,\,0,\,0,\,0,\,0,\,0,\,0,\,0,\,0,\,0,\,0,\,1,\,0,\,0,\,0,\,0,\,0\right) ,x \rangle \leq \frac{2361185035008561037417}{2361183241434822606848} \\
-\frac{4722372222305608191517}{4722366482869645213696} &\leq \langle \left(0,\,0,\,0,\,0,\,0,\,0,\,0,\,0,\,0,\,0,\,0,\,0,\,0,\,0,\,0,\,0,\,0,\,0,\,1,\,0\right) ,x \rangle \leq \frac{4722372222305608191517}{4722366482869645213696} \\
-\frac{2361185124158778945083}{2361183241434822606848} &\leq \langle \left(0,\,0,\,0,\,0,\,0,\,0,\,1,\,0,\,0,\,0,\,0,\,0,\,0,\,0,\,0,\,0,\,0,\,0,\,0,\,0\right) ,x \rangle \leq \frac{2361185124159496218833}{2361183241434822606848} \\
\end{align*}}
\subsection{Counterexample based on inequality $32$}
If $H$ is the polyhedron obtained from $N$ by removing inequality $I_{32}$, then $H$ is a bounded Hirsch polytope and $N = H \cap C$, where $C$ is the combinatorial $20$-cube defined by
{\footnotesize\begin{align*}
-\frac{590297603932444082281}{590295810358705651712} &\leq \langle \left(0,\,0,\,0,\,0,\,0,\,0,\,0,\,0,\,0,\,0,\,0,\,0,\,0,\,0,\,0,\,0,\,1,\,0,\,0,\,-\frac{1}{100}\right) ,x \rangle \leq \frac{1180595210734606146051}{1180591620717411303424} \\
-\frac{147579558414782802697}{147573952589676412928} &\leq \langle \left(0,\,0,\,0,\,0,\,0,\,0,\,0,\,0,\,0,\,0,\,0,\,0,\,0,\,1,\,0,\,0,\,0,\,0,\,0,\,-\frac{1}{10000}\right) ,x \rangle \leq \frac{147579558414782802697}{147573952589676412928} \\
-\frac{2361185124158778945083}{2361183241434822606848} &\leq \langle \left(0,\,0,\,0,\,0,\,0,\,0,\,1,\,0,\,0,\,0,\,0,\,0,\,0,\,0,\,0,\,0,\,0,\,0,\,0,\,0\right) ,x \rangle \leq \frac{2361185124159496218833}{2361183241434822606848} \\
-\frac{9444738705175253405213}{9444732965739290427392} &\leq \langle \left(0,\,0,\,0,\,0,\,0,\,0,\,0,\,0,\,0,\,0,\,0,\,0,\,0,\,0,\,0,\,0,\,0,\,1,\,0,\,-\frac{1}{10}\right) ,x \rangle \leq \frac{4722369467376345962163}{4722366482869645213696} \\
-\frac{2417851643085731713983019}{2417851639229258349412352} &\leq \langle \left(0,\,0,\,0,\,0,\,0,\,0,\,0,\,0,\,1,\,0,\,0,\,0,\,0,\,0,\,0,\,0,\,0,\,0,\,0,\,0\right) ,x \rangle \leq \frac{2417851643085731713983019}{2417851639229258349412352} \\
-\frac{40066579960122721}{36028797018963968} &\leq \langle \left(0,\,0,\,0,\,1,\,0,\,0,\,0,\,0,\,0,\,0,\,0,\,0,\,0,\,0,\,0,\,0,\,0,\,0,\,0,\,-\frac{1}{500000000000}\right) ,x \rangle \leq \frac{20040505775302443}{18014398509481984} \\
-\frac{4276905141513425}{8796093022208} &\leq \langle \bigg(100,\,27,\,0,\,-\frac{1}{500},\,\frac{1}{80},\,0,\,0,\,0,\,0,\,0,\,0,\,0,\,0,\,-100000,\,-10000000,\,-10000000,\\&\,-10000000,\,-100000000,\,-100000000,\,-1000000000\bigg) ,x \rangle \leq 1 \\
-\frac{2417851645399615732725419}{2417851639229258349412352} &\leq \langle \left(0,\,0,\,0,\,0,\,0,\,0,\,0,\,0,\,0,\,0,\,0,\,0,\,1,\,0,\,0,\,0,\,0,\,0,\,0,\,0\right) ,x \rangle \leq \frac{2417851645399615732725419}{2417851639229258349412352} \\
-\frac{39614081257140112583132697033}{39614081257132168796771975168} &\leq \langle \left(0,\,0,\,0,\,0,\,0,\,0,\,0,\,0,\,0,\,0,\,1,\,0,\,0,\,0,\,0,\,0,\,0,\,0,\,0,\,0\right) ,x \rangle \leq \frac{9903520314285028145783175131}{9903520314283042199192993792} \\
-\frac{76859732694334399}{72057594037927936} &\leq \langle \left(0,\,1,\,0,\,0,\,0,\,0,\,0,\,0,\,0,\,0,\,0,\,0,\,0,\,0,\,0,\,0,\,0,\,0,\,0,\,\frac{27}{1000000000}\right) ,x \rangle \leq \frac{77296343773163779}{72057594037927936} \\
-\frac{291112682551516277}{288230376151711744} &\leq \langle \left(1,\,0,\,0,\,0,\,0,\,0,\,0,\,0,\,0,\,0,\,0,\,0,\,0,\,0,\,0,\,0,\,0,\,0,\,0,\,\frac{1}{10000000}\right) ,x \rangle \leq \frac{582225359826457723}{576460752303423488} \\
-\frac{4722372222305608191517}{4722366482869645213696} &\leq \langle \left(0,\,0,\,0,\,0,\,0,\,0,\,0,\,0,\,0,\,0,\,0,\,0,\,0,\,0,\,0,\,0,\,0,\,0,\,1,\,-\frac{1}{10}\right) ,x \rangle \leq \frac{4722372337094327451073}{4722366482869645213696} \\
-\frac{9671406563106733594076309}{9671406556917033397649408} &\leq \langle \left(0,\,0,\,0,\,0,\,0,\,0,\,0,\,0,\,0,\,0,\,0,\,1,\,0,\,0,\,0,\,0,\,0,\,0,\,0,\,0\right) ,x \rangle \leq \frac{9671406563106733594076309}{9671406556917033397649408} \\
-\frac{309485009827534768921207957}{309485009821345068724781056} &\leq \langle \left(0,\,0,\,0,\,0,\,0,\,0,\,0,\,0,\,0,\,1,\,0,\,0,\,0,\,0,\,0,\,0,\,0,\,0,\,0,\,0\right) ,x \rangle \leq \frac{309485009827534768921207957}{309485009821345068724781056} \\
-\frac{1180593414291149733993}{1180591620717411303424} &\leq \langle \left(0,\,0,\,0,\,0,\,0,\,0,\,0,\,0,\,0,\,0,\,0,\,0,\,0,\,0,\,0,\,1,\,0,\,0,\,0,\,-\frac{1}{100}\right) ,x \rangle \leq \frac{590296708580433857741}{590295810358705651712} \\
-\frac{1180593503442084915409}{1180591620717411303424} &\leq \langle \left(0,\,0,\,0,\,0,\,0,\,0,\,0,\,1,\,0,\,0,\,0,\,0,\,0,\,0,\,0,\,0,\,0,\,0,\,0,\,0\right) ,x \rangle \leq \frac{1180593503442084915409}{1180591620717411303424} \\
-\frac{151699951751958945}{144115188075855872} &\leq \langle \left(0,\,0,\,1,\,0,\,0,\,0,\,0,\,0,\,0,\,0,\,0,\,0,\,0,\,0,\,0,\,0,\,0,\,0,\,0,\,0\right) ,x \rangle \leq \frac{80781273703537389}{72057594037927936} \\
-\frac{2392514266058385}{2251799813685248} &\leq \langle \left(0,\,0,\,0,\,0,\,1,\,0,\,0,\,0,\,0,\,0,\,0,\,0,\,0,\,0,\,0,\,0,\,0,\,0,\,0,\,\frac{1}{80000000000}\right) ,x \rangle \leq \frac{76542601517243887}{72057594037927936} \\
-\frac{2361185035008561037417}{2361183241434822606848} &\leq \langle \left(0,\,0,\,0,\,0,\,0,\,0,\,0,\,0,\,0,\,0,\,0,\,0,\,0,\,0,\,1,\,0,\,0,\,0,\,0,\,-\frac{1}{100}\right) ,x \rangle \leq \frac{9444740162991988001579}{9444732965739290427392} \\
-\frac{18889473462374272866723}{18889465931478580854784} &\leq \langle \left(0,\,0,\,0,\,0,\,0,\,1,\,0,\,0,\,0,\,0,\,0,\,0,\,0,\,0,\,0,\,0,\,0,\,0,\,0,\,0\right) ,x \rangle \leq \frac{18889473462375840754951}{18889465931478580854784} 
\end{align*}}

\vskip16pt
\noindent K.~P.~Fallon, M.~Janusiak, E.~D.~Kim*, A.~McLain

\vskip16pt
\noindent *Corresponding author: {\tt ekim@uwlax.edu}, 1725 State Street, La Crosse, WI 54601 (USA)

\vskip16pt
\noindent Keywords: cube, linear programming, polyhedron

\end{document}